\documentclass{article}
\usepackage{amssymb,amsmath}
\usepackage[mathscr]{euscript}
\input{epsf}

\newcounter{sec}

\newcounter{punct}[sec]

\def\punct{\refstepcounter{punct}{\arabic{sec}.\arabic{punct}.  }}

\newtheorem{theorem}{Theorem}[sec]
\newtheorem{proposition}[theorem]{Proposition}

\newtheorem{lemma}[theorem]{Lemma}

\newtheorem{corollary}[theorem]{Corollary}

\newtheorem{question}{Question}[sec]

\def\COUNTERS{\addtocounter{sec}{1}
              \setcounter{punct}{0}
          \setcounter{equation}{0}
          \setcounter{theorem}{0}
          }
          
          \def\sm{\smallskip}

\begin{document}

\newcommand{\supp}{\mathop {\mathrm {supp}}\nolimits}
\newcommand{\rk}{\mathop {\mathrm {rk}}\nolimits}
\newcommand{\Aut}{\mathop {\mathrm {Aut}}\nolimits}
\newcommand{\Out}{\mathop {\mathrm {Out}}\nolimits}
\renewcommand{\Re}{\mathop {\mathrm {Re}}\nolimits}

\def\0{\mathbf 0}

\def\ov{\overline}
\def\wh{\widehat}
\def\wt{\widetilde}

\renewcommand{\rk}{\mathop {\mathrm {rk}}\nolimits}
\renewcommand{\Aut}{\mathop {\mathrm {Aut}}\nolimits}
\renewcommand{\Re}{\mathop {\mathrm {Re}}\nolimits}
\renewcommand{\Im}{\mathop {\mathrm {Im}}\nolimits}
\newcommand{\sgn}{\mathop {\mathrm {sgn}}\nolimits}

\def\bfa{\mathbf a}
\def\bfb{\mathbf b}
\def\bfc{\mathbf c}
\def\bfd{\mathbf d}
\def\bfe{\mathbf e}
\def\bff{\mathbf f}
\def\bfg{\mathbf g}
\def\bfh{\mathbf h}
\def\bfi{\mathbf i}
\def\bfj{\mathbf j}
\def\bfk{\mathbf k}
\def\bfl{\mathbf l}
\def\bfm{\mathbf m}
\def\bfn{\mathbf n}
\def\bfo{\mathbf o}
\def\bfp{\mathbf p}
\def\bfq{\mathbf q}
\def\bfr{\mathbf r}
\def\bfs{\mathbf s}
\def\bft{\mathbf t}
\def\bfu{\mathbf u}
\def\bfv{\mathbf v}
\def\bfw{\mathbf w}
\def\bfx{\mathbf x}
\def\bfy{\mathbf y}
\def\bfz{\mathbf z}

\def\bfA{\mathbf A}
\def\bfB{\mathbf B}
\def\bfC{\mathbf C}
\def\bfD{\mathbf D}
\def\bfE{\mathbf E}
\def\bfF{\mathbf F}
\def\bfG{\mathbf G}
\def\bfH{\mathbf H}
\def\bfI{\mathbf I}
\def\bfJ{\mathbf J}
\def\bfK{\mathbf K}
\def\bfL{\mathbf L}
\def\bfM{\mathbf M}
\def\bfN{\mathbf N}
\def\bfO{\mathbf O}
\def\bfP{\mathbf P}
\def\bfQ{\mathbf Q}
\def\bfR{\mathbf R}
\def\bfS{\mathbf S}
\def\bfT{\mathbf T}
\def\bfU{\mathbf U}
\def\bfV{\mathbf V}
\def\bfW{\mathbf W}
\def\bfX{\mathbf X}
\def\bfY{\mathbf Y}
\def\bfZ{\mathbf Z}

\def\frD{\mathfrak D}
\def\frL{\mathfrak L}
\def\frG{\mathfrak G}
\def\frg{\mathfrak g}
\def\frh{\mathfrak h}
\def\frf{\mathfrak f}
\def\frl{\mathfrak l}

\def\bfw{\mathbf w}

\def\R {{\mathbb R }}
 \def\C {{\mathbb C }}
  \def\Z{{\mathbb Z}}
  \def\H{{\mathbb H}}
\def\K{{\mathbb K}}
\def\N{{\mathbb N}}
\def\Q{{\mathbb Q}}
\def\A{{\mathbb A}}
\def\O {{\mathbb O }}

\def\T{\mathbb T}
\def\P{\mathbb P}

\def\G{\mathbb G}

\def\cD{\EuScript D}
\def\cL{\mathscr L}
\def\cK{\EuScript K}
\def\cM{\EuScript M}
\def\cN{\EuScript N}
\def\cP{\EuScript P}
\def\cQ{\EuScript Q}
\def\cR{\EuScript R}
\def\cW{\EuScript W}
\def\cY{\EuScript Y}
\def\cF{\EuScript F}
\def\cG{\EuScript G}

\def\bbA{\mathbb A}
\def\bbB{\mathbb B}
\def\bbD{\mathbb D}
\def\bbE{\mathbb E}
\def\bbF{\mathbb F}
\def\bbG{\mathbb G}
\def\bbI{\mathbb I}
\def\bbJ{\mathbb J}
\def\bbL{\mathbb L}
\def\bbM{\mathbb M}
\def\bbN{\mathbb N}
\def\bbO{\mathbb O}
\def\bbP{\mathbb P}
\def\bbQ{\mathbb Q}
\def\bbS{\mathbb S}
\def\bbT{\mathbb T}
\def\bbU{\mathbb U}
\def\bbV{\mathbb V}
\def\bbW{\mathbb W}
\def\bbX{\mathbb X}
\def\bbY{\mathbb Y}

\def\kappa{\varkappa}
\def\epsilon{\varepsilon}
\def\phi{\varphi}
\def\le{\leqslant}
\def\ge{\geqslant}

\def\B{\mathrm B}

\def\la{\langle}
\def\ra{\rangle}
\def\tri{\triangleright}

\def\lambdA{{\boldsymbol{\lambda}}}
\def\alphA{{\boldsymbol{\alpha}}}
\def\betA{{\boldsymbol{\beta}}}
\def\mU{{\boldsymbol{\mu}}}

\def\const{\mathrm{const}}
\def\rem{\mathrm{rem}}
\def\even{\mathrm{even}}
\def\SO{\mathrm{SO}}
\def\SL{\mathrm{SL}}
\def\GL{\operatorname{GL}}
\def\End{\operatorname{End}}
\def\Mor{\operatorname{Mor}}
\def\Aut{\operatorname{Aut}}
\def\inv{\operatorname{inv}}
\def\red{\operatorname{red}}
\def\Ind{\operatorname{Ind}}
\def\dom{\operatorname{dom}}
\def\im{\operatorname{im}}
\def\md{\operatorname{mod\,}}

\def\ZZ{\mathbb{Z}_{p^\mu}}
\def\F{\mathbb{F}}

\def\cH{\EuScript{H}}
\def\cQ{\EuScript{Q}}
\def\cL{\EuScript{L}}
\def\cX{\EuScript{X}}

\def\Di{\Diamond}
\def\di{\diamond}

\def\fin{\mathrm{fin}}
\def\ThetA{\boldsymbol {\Theta}}

\def\0{\boldsymbol{0}}

\begin{center}

{\bf\Large Description of unitary representations

\medskip

of the group of infinite $p$-adic integer matrices}

\bigskip

\sc \large
 Yury A. Neretin%
\footnote{Supported by the grants FWF, P28421, P31591.}

\end{center}

{\small We classify  irreducible unitary representations of the group of all infinite matrices over a $p$-adic field ($p\ne 2$) with integer elements
 equipped with a  natural topology. Any irreducible representation passes
through a group $\GL$ of infinite matrices over a residue ring modulo $p^k$.
Irreducible representations of the latter group are induced from
finite-dimensional representations of certain open subgroups.}

\section{Introduction}

\COUNTERS

{\bf\punct Notations and definitions.}
a) {\sc Rings.}
Let $p$ be a prime,
$$
p>2.
$$
 Let $\Z_{p^n} := \Z/p^n\Z$ be a residue
ring, $\F_p := \Z_p$ be the field with $p$ elements. 
The ring of $p$-adic integers $\O_p$ is
the projective limit
$$
\O_p=
\lim\limits_{\longleftarrow n} \Z_{p^n}
$$
of the following chain (see, e.g.,
\cite{Ser}:
$$\dots \longleftarrow 
\Z_{p^{n-1}} \longleftarrow \Z_{p^n} \longleftarrow \Z_{p^{n+1}} \longleftarrow \dots
$$
we have $\Z_{p^n}=\O_p/p^n\O_p$. Denote by $\Q_p$ the field of $p$-adic numbers.

\sm 

b) {\sc The infinite symmetric group and oligomorphic groups.}
Let $\Omega$ be a countable set. Denote by $S(\Omega)$ the
group of all permutations of $\Omega$, denote $S_\infty:=S(\N)$. The topology on
the {\it infinite symmetric group} $S(\Omega)$
 is determined by the
condition: stabilizers of finite subsets are open subgroups and these subgroups
form a fundamental system of neighborhoods of the unit%
\footnote{Thus we get a structure of a Polish group, moreover
	 this topology is a unique separable topology on the
	infinite symmetric group, see
	\cite{KR}. In particular, this means that a unitary representation of $S_\infty$
	in a separable Hilbert space is automatically continuous.}.
 Equivalently, a sequence $g^{(\alpha)}$ converges to $g$ if for each 
 $\omega \in \Omega$ we have $\omega g^{(\alpha)} = \omega g$ for sufficiently
large $\alpha$.

A closed subgroup $G$ of $S(\Omega)$ is called {\it oligomorphic} if for each $k$ it has only a
finite number of orbits on the product $\Omega \times\dots\times \Omega$ of $k$ copies of $\Omega$, see
\cite{Cam}.

\sm

c) {\sc Modules $\frl(\Z_{p^n})$ and groups $\GL(\infty, \Z_{p^n})$.}
Define the module $\frl(\Z_{p^n})$ as the set of all sequences 
$v = (v_1, v_2,\dots)$, where $v_j \in \Z_{p^n}$
and $v_j = 0$ for sufficiently large $j$. The set $\frl(\Z_{p^n})$ 
is countable, we
equip it
with a discrete topology. Denote by $e_j$ the standard basis elements,
i.e., $e_j$ has a unit on $j$-th place, other elements are $0$.

Define groups $\GL(\infty, \Z_{p^n})$
 as groups of infinite invertible matrices $g$
over $\Z_{p^n}$ such that:

\sm 

$\bullet$ each row of $g$ contains only a
finite number of nonzero elements;

\sm

$\bullet$  each column contains only a
finite number of nonzero elements;

\sm

$\bullet$ the inverse matrix $g^{-1}$ satisfies the same conditions.

\sm 

Notice that rows of a matrix $g$ are precisely vectors $e_i g$, and columns are
$e_j g^t$ (we denote by $g^t$ a {\it transposed matrix}).

\sm 

 {\it Actually, the topic of this paper are representations of $\GL(\infty, \Z_{p^n})$.}
 
 \sm
 
 This group is continual and we must define a {\it topology} on $\GL(\infty, \Z_{p^n})$.
 A
sequence $g^{(\alpha)}\in   \GL(\infty, \Z_{p^n})$ converges to $g$ if all sequences $e_ig^{(\alpha)}$ and $e_i(g^{(\alpha)})^t$
are eventually constant and their limits are $e_ig$ and $e_jg^t$ respectively.
Thus we get a structure of a totally disconnected topological group.

The group $\GL(\infty, \Z_{p^n})$ acts on the countable set $\frl(\Z_{p^n})\oplus \frl(\Z_{p^n})$ by transformations
$$
(v,w)\mapsto (vg, wg^{t-1}).
$$
In particular, this define an embedding of $\GL(\infty, \Z_{p^n})$  to a symmetric group $S\bigl(\frl(\Z_{p^n})\oplus \frl(\Z_{p^n})\bigr)$.
The image of the group $\GL(\infty, \Z_{p^n})$ is a closed subgroup of $S\bigl(\frl(\Z_{p^n})\oplus \frl(\Z_{p^n})\bigr)$ and the induced
topology coincides with the natural topology on $\GL(\infty, \Z_{p^n})$. By
\cite{Ner-p}, Lemma 3.7, {\it the
group $\GL(\infty, \Z_{p^n})$  is oligomorphic}.


\sm

d) {\sc Modules $\frl(\O_p)$ and groups $\GL(\infty, \O_p)$}. Denote by $\frl (\O_p)$ the set of all sequences 
$r = (r_1, r_2, \dots )$,
where $r_j \in  \O_p$ and $|r_j| \to 0$ as $j \to \infty$.
 The space $\frl(\O_p)$ is a projective limit,
$$
\frl(\O_p)=\lim\limits_{\longleftarrow n}\frl(\Z_{p^n}),
$$
we equip it with the {\it topology of the projective limit}. In other words, a sequence
$r^{(j)}\in \frl(\O_p)$ converges if for any $p^n$ the reduction of $r^{(j)}$ modulo $p^n$ is eventually
constant in $\Z_{p^n}$.

\sm

 We define $\GL(\infty, \O_p)$
 as the group of all
infinite matrices $g$ over $\O_p$ such that:

\sm

$\bullet$ each row of $g$ is an element of $\frl(\O_p)$;

\sm

$\bullet$ each column of $g$ is an element of $\frl(\O_p)$;

\sm

$\bullet$ the matrix $g$ has an inverse and $g^{-1}$ satisfies the same conditions.

\sm

 We say that a sequence $g^{(\alpha)} \in \GL(\infty, \O_p)$
 converges to $g$ if for any $i$ the
sequence $e_i g^{(\alpha)}$ converges to $e_i g$ and for any $j$
 the sequence $e_i (g^{(\alpha)})^t$ converges
to $e_j g^t$.
This determines a structure of a totally disconnected topological group on
$\GL(\infty, \O_p)$.

We have obvious homomorphisms $\GL(\infty, \Z_{p^n}) \to \GL(\infty, \Z_{p^{n-1}})$, the group
$\GL(\infty, \O_p)$ is the projective limit
$$
\GL(\infty, \O_p) = \lim_{\longleftarrow n}
\GL(\infty, \Z_{p^n})$$
and its topology is the topology of projective limit.

\sm

{\bf\punct Preliminary remarks.} A priori we know the following statement:

\begin{theorem}
	\label{th:oligo} {\rm a)}
	 The group $\GL(\infty, \O_p)$
	  is a type $I$ group, it has a countable number of irreducible unitary representations. Any unitary representation
	$\GL(\infty, \O_p)$ is a sum of irreducible representations. Any irreducible unitary representation of $\GL(\infty, \O_p)$ 
	is in fact a representation of some group $\GL(\infty, \Z_{p^n})$.
	
	\sm
	
	{\rm b)} Each irreducible representation of $\GL(\infty, \Z_{p^n})$ 
	is induced from a finite-dimensional representation of an open subgroup. More precisely, for any irreducible unitary representation of $\GL(\infty, \Z_{p^n})$ there exists an open subgroup
	$\wh Q\subset \GL(\infty, \Z_{p^n})$, a normal subgroup $Q\subset\wh Q$
	of finite index and an irreducible
	representation $\nu$ of $\wh Q$, which is trivial on $Q$, such that $\rho$ is induced from $\nu$.
\end{theorem}

This is a special case of a theorem of Tsankov about unitary representations of oligomorphic groups
and projective limits of holomorphic groups, see 
\cite{Tsa}, Theorem 1.3%
\footnote{A reduction of representations of $\GL(\infty,\O_p)$
to representations of quotients $\GL(\infty,\ZZ)$
easily follows from   \cite{Ner-book}, Proposition VII.1.3,
see \cite{Ner-p}, Corollary 3.5. In our proof of Theorem \ref{th:main}
Tsankov's theorem  is used in the proof of Proposition \ref{pr:theta},
which was done in \cite{Ner-p}.}%
. It seems that
\cite{Tsa},
\cite{BT}
 is not sufficient to give a precise answer in our case.
 
 \sm 
 
 Let us give a definition of an {\it induced representation}
 (see, e.g., \cite{Ser-rep}, Sect. 7, Kirillov, Sect. 13), which is appropriate in
 our case. Let $G$ be a totally disconnected separable group, $Q$ its open subgroup.
 Let $\nu$ be a unitary representation of $Q$ in a Hilbert space $V$. Consider the space
 $H$ of $V$-valued functions $f$ on a countable homogeneous space $Q\setminus G$ such that 
 $$\sum_{x\in Q\setminus G} \|f(x)\|^2 < \infty.$$
 Equip $H$ with the
  inner product
 $$
 \la f_1, f_2\ra_H := 
 \sum_{x\in Q\setminus G}
 \la f_1(x), f_2(x)\ra_V.
 $$
 Let $U$ be a function on $G\times (Q\setminus G)$ taking values in the group of unitary operators
 in $V$ such that:
 
 \sm
 
 $\bullet$ Formula 
 $$\rho(g)f(x) = U(g, x)f(xg)
 $$
  determines a representation of $G$ in
 $H$.
 
 \sm
 
  $\bullet$ Let $x_0$ be the initial point of $Q\setminus G$, i.e., $x_0 Q = x_0$.
   Then for $q \in  Q$ we
 have $U(q, x_0) = \nu(q)$.
 
 \sm
 
 The first condition implies that the function $U(g,x)$ 
 satisfies the functional equation
 $$
 U(x,g_1g_2)=U(x,g_1)\, U(xg_1,g_2).
 $$
 It can be shown that $U(g,x)$
  is uniquely defined up to a natural calibration
 $$U(g, x) \sim A(gx)^{-1}U(g, x)A(x),$$
 where  $A$ is a function on $Q\setminus G$ taking values in the unitary group of $V$ (see, e.g.,
 \cite{Kir}, Sect 13.1). For this reason, an induced representation
 $\rho(g) = \Ind^G_Q(\nu)$
 is canonically defined up to a unitary equivalence.
 
 We also can choose $U(g, x)$ in the following way. For any $x \in Q\setminus G$ we choose
 an element $s(x)\in G$ such that $x_0 s(x) = x$. Then $U(g, x) = \nu(q)$, where $q$ is
 determined from the condition $s(x)g = q\,s(xg)$.

\sm

{\bf \punct The statement.%
\label{ss:statement}} Thus we fix a ring $\Z_{p^\mu}$ and examine the group
$$\G := \GL(\infty, \Z_{p^\mu}).$$
We consider two right actions of $\G$ on $\frl(\ZZ), g : v\mapsto vg$,
$ g :v\mapsto vg^{t-1}$. Define a
pairing
$$
\frl(\ZZ)\times \frl(\ZZ)\to\ZZ
$$
by
\begin{equation}
\{v, w\} := \sum v_j w_j= v w^t,
\label{eq:pairing}
\end{equation}
our action preserves this pairing, i.e.,
$$\{ vg, vg^{t-1}\} = \{v, w\}.$$

Let $L \subset\frl(\ZZ)$, $M \subset'\frl(\ZZ)$ be finitely generated 
$\ZZ$-submodules. Denote
by $\wh \G[L; M]$ the subgroup of $\G$ consisting of $g$ sending $L$ to 
$L$ and $M$ to $M$.
By $\G^\circ[L; M] \subset  \wh \G[L; M]$
we denote group of matrices fixing $L$ and $M$ pointwise.
Obviously, the quotient group $ \wh \G[L; M]/\G^\circ[L; M]$ is finite, it acts on
the direct sum $L\oplus M$ preserving the pairing $\{f,g\}$. Any irreducible representation
$\tau$ of $ \wh \G[L; M]/\G^\circ[L; M]$
 can be regarded as a representation $\wh \tau$
 the group 
$\wh \G[L; M]$, which
is trivial on $\G^\circ[L; M]$. For given $L$, $M$, $\tau$ we consider the representation
$$
\Ind_{\wh \G[L; M]}^\G(\wh\tau)
$$
of $\G$ induced from the representation $\wh \tau$ of the group $\wh \G[L; M]$.
Olshanski \cite{Olsh} obtained the following statement%
\footnote{A proof in \cite{Olsh} is only sketched, other proofs were given by Dudko \cite{Dud}
and Tsankov \cite{Tsa}.}
for the group
$\GL(\infty,\F_p)=\GL(\infty,\Z_p)$.

\begin{theorem}
\label{th:main-o}  {\rm a)} Any irreducible unitary representation of the group 
$\GL(\infty,\F_p)$
has this form.

\sm

{\rm b)} Two irreducible representations can be equivalent only for a trivial reason,
i.e.,
$$
\Ind_{\wh \G[L_1; M_1]}^\G(\tau_1)\sim
\Ind_{\wh \G[L_2; M_2]}^\G(\tau_2)
$$
if and only if there exists $h \in \G$ such that $L_1h = L_2h$, 
$M_1h^{t-1} =M_2 h^{t-1}$ and
$\tau_2(q) = \tau_1(hqh^{-1})$.	
\end{theorem}

For groups $\GL(\infty,\ZZ)$ with $\mu>1$  the situation is more delicate.
Let $L$, $M$ actually be contained in $(\ZZ)^m\subset \frl (\ZZ)$.
Fix a matrix $b$ such that%
\footnote{We assume that each row of $b$ and each column of $c$
contains only a finite number of nonzero elements.} $\ker b=L$ and a matrix $c$ such that
$\ker c^t=M$.

\begin{lemma}
	\label{l:G-circ}
The group $\G^\circ[L;M]$ consists of all invertible 
 matrices admitting the following
representation as a block matrix of size $m+\infty$:
\begin{equation}
\label{eq:G-circ}
g=\begin{pmatrix}
a&bv\\w c&z
\end{pmatrix},
\end{equation}
where the block '$a$' can be written in   both forms
$$
a=1-b S,\qquad a=1-T c.
$$
\end{lemma}

Next, define a subgroup $\G^\bullet[L;M]\subset \G^\circ[L;M]$
consisting of matrices having the form
\begin{equation}
\label{eq:G-bullet}
g=\begin{pmatrix}
1-buc&bv\\w c&z
\end{pmatrix}.
\end{equation}

\begin{proposition}
	\label{pr:bullet}
The group $\G^\bullet[L;M]$ is the minimal
subgroup of finite index in $\wh\G[L;M]$,
i.e., it is contained in any subgroup of finite index in $\wh\G[L;M]$.
\end{proposition}

\begin{theorem}
	\label{th:main}
	{\rm a)}  Any irreducible unitary representation of
	$\G$ is induced from a representation $\tau$ of some group
	$\wh\G[L;M]$ that is trivial on the subgroup $\G^\bullet[L;M]$.
	
	\sm
	
	{\rm b)} Two irreducible  representations of this kind can be
	 equivalent only for the trivial reason as in Theorem {\rm \ref{th:main-o}}.
\end{theorem}

{\sc Remark.} Recall that $p\ne 2$. In several places of our proof we divide elements
of residue rings $\ZZ$ by 2. Usually, this division can be replaced by 
longer considerations. But  in Lemma \ref{l:2} this seems crucial.
\hfill $\boxtimes$

\sm

{\sc Remark.}
Let $L$, $M\subset p\cdot\frl(\ZZ)$. Then
$\G[L;M]$ contains a congruence subgroup
$N$ consisting of elements of $\G$ that are equal $1$ modulo $p^{\mu-1}$.
Since $N$ is a normal subgroup
in $\G$, it is normal in $\wh\G[L;M]$.
Let $\tau$ be trivial on $N$.
Then the
induced representation $\Ind^\G_{\G[L;M]}(\wh\tau)$ is trivial
on the congruence subgroup $N$
and actually we get representations of $\GL(\infty,\Z_{p^{\mu-1}})$.
\hfill $\boxtimes$

\sm

{\sc Remark.} The statement b) is a general fact for oligomorphic groups,
see \cite{Tsa}, Proposition 4.1(ii). So we omit a proof (in our case this
can be easily established by examination of intertwining operators).
\hfill $\boxtimes$

\sm

{\bf\punct  Remarks. Infinite-dimensional $p$-adic groups}. Now there exists
a well-developed representation theory of infinite symmetric groups 
and of infinite-dimensional real classical groups.
Parallel development in the $p$-adic case meet some difficulties. However,
infinite dimensional
 $p$-adic
group were  a topic of sporadic attacks
since late 1980s, see
\cite{NNO},
\cite{Zel},
\cite{Naz}. We indicate some works on $p$-adic groups
and their parallels with
 nontrivial constructions for real and symmetric groups.

\sm

a) An extension of the Weil representation of the infinite-dimensional symplectic group $\mathrm{Sp}(2\infty,\C)$
to the semigroup of lattices (Nazarov
\cite{NNO},
\cite{Naz}, see a partial
exposition in
\cite{Ner-gauss}, Sect. 11.1-11.2).

\sm

b) A construction of projective limits of $p$-adic Grassmannians and quasiinvariant actions of $p$-adic $\GL(\infty)$ on these Grassmannians
\cite{Ner-hua-p}. This is an
analog of virtual permutations (or Chinese restaurant process, see, e.g., \cite{Ald}, 11.19, they are a base of
harmonic analysis related to infinite symmetric group, see
\cite{KOV}), and of projective limits of compact symmetric spaces (see
\cite{Pick},
\cite{Ner-hua}), they are
a standpoint for a harmonic analysis related to infinite-dimensional classical
groups, see
\cite{BO}).

\sm

c) An attempt to describe a multiplication of double cosets (see the next
section) for $p$-adic classical groups in
\cite{Ner-build}. In any case this leads to a strange geometric construction, namely to simplicial maps of Bruhat--Tits buildings whose
	boundary values  are rational maps
	of $p$-adic Grassmannians.
	
	\sm
	
	d) The work
	\cite{Buf}
	 contain a $p$-adic construction in the spirit of exchangeability%
	\footnote{i.e., of higher analogs of the de Finetti theorem, see,
	\cite{Ald}}, namely, descriptions of
	invariant ergodic measures on spaces of infinite $p$-adic matrices. 
	By the Wigner--Mackey trick (see, e.g., \cite{Kir}, Sect. 13.3), such kind of
	statements can be translated to a description of spherical functions on certain
	groups.
	
	\sm
	
	So during last years new elements of a nontrivial picture
	related to
	infinite-dimensional $p$-adic groups appeared. For this reason,
	 understanding of representations $\GL(\infty,\O_p)$ becomes necessary.
	 
	 \sm
	 
{\bf \punct Another completion of a group of infinite matrices over $\Z_{p^n}$.}
Define a group $\cG$ consisting of infinite matrices $g$ over $\Z_{p^n}$
such that: 

\sm

$\bullet$ $g$ contains
only a finite number of elements in each column;

\sm

$\bullet$ $g^{-1}$ exists and satisfies the same property.

\sm

A sequence $g^{(\alpha)}$ converges to $g$ if for each $j$ we have a convergence of $e_j g^{(\alpha)}$. 

\sm

Clearly, $\cG\supset \G$. Classification of irreducible unitary representations of $\cG$ is the following.
For each finitely generated submodule in $\frl(\Z_{p^n})$ we consider the subgroup $\wh\cG[L]$ consisting
of transformations sending $L$ to itself and the subgroup $\cG^\circ[L]$ fixing $L$ pointwise.

\begin{proposition}
 Any irreducible unitary representation of $\cG$ is induced from a representation of some group
 $\wh\cG[L]$ trivial on $\cG^\circ[L]$.
\end{proposition}

This follows from Theorem \ref{th:main}, on the other hand this can be deduced in a straightforward way from
Tsankov's result \cite{Tsa}.

	 \section{Preliminaries: the category of double cosets}
	 
	 \COUNTERS

{\bf\punct Multiplication of double cosets and the category $\cK$.%
\label{ss:cosets}} Here we
discuss a version of a general construction of multiplication of double cosets (see \cite{Olsh-howe},
\cite{Olsh},
\cite{Ner-book},
\cite{Ner-symm},
\cite{Ner-p}).

Denote by $\G_\fin \subset \G$ the subgroup of {\it finitary} matrices, i.e., matrices $g$ such
that $g-1$ has only a finite number of nonzero elements. For $\alpha = 0$, 1, 
\dots denote by
$\G(\alpha) \subset \G$ the subgroups consisting of matrices having the form 
$\begin{pmatrix}
	1_\alpha&0\\0&u
\end{pmatrix}$, where
$1_\alpha$ denotes the unit matrix of  size $\alpha$ and $u$ is an arbitrary invertible matrix
over $\ZZ$. Obviously, $\G(\alpha)$ is isomorphic to 
$\G$. Consider double cosets spaces
$\G(\alpha) \setminus \G/\G(\beta)$, their elements are matrices determined up to the equivalence
\begin{equation}
\begin{pmatrix}
a&b\\c&d
\end{pmatrix}\sim
\begin{pmatrix}
1_\alpha&0\\0&u
\end{pmatrix}
\begin{pmatrix}
a&b\\c&d
\end{pmatrix}
\begin{pmatrix}
1_\beta&0\\0&v
\end{pmatrix}
=
\begin{pmatrix}
a&bv\\uc&udv
\end{pmatrix},
\label{eq:sim}
\end{equation}
where a matrix $g$ is represented as a block matrix of  size $(\alpha+\infty)\times(\beta +\infty)$.
For a matrix $g$ we write the corresponding double coset as
$$
\left[
\begin{array}{c|c}
a&b\\\hline c&d
\end{array}
\right]_{\alpha\beta},
$$
we will omit subscripts $\alpha\beta$ if it is not necessary to indicate a size.
We wish to define a natural multiplication
$$\G(\alpha) \setminus \G/\G(\beta)\,\, \times\,\,
 \G(\beta) \setminus \G/\G(\gamma)\, \to\, \G(\alpha) \setminus \G/\G(\gamma).
$$
Let $\frg_1 \in \G(\alpha)\setminus\G/\G(\beta)$, 
$\frg_2 \in \G(\beta)\setminus \G/\G(\gamma)$ be double cosets. By
\cite{Ner-p}, Lemma
4.1, any double coset has a representative in $\G_\fin$. Choose such representatives
$g_1$ and $g_2$ for $\frg_1$, $\frg_2$,
\begin{equation}
g_1=\left[\begin{array}{c|cc}
a&b&\\
\hline
c&d&\\
&&1_\infty
\end{array}\right]_{\alpha\beta},
\,\,
g_1=\left[\begin{array}{c|cc}
p&q&\\
\hline
r&t&\\
&&1_\infty
\end{array}\right]_{\beta\gamma}
.
\label{eq:g1g2}
\end{equation}
let sizes of submatrices $\begin{pmatrix}
a&b\\c&d
\end{pmatrix}$, $\begin{pmatrix}
	p&q\\r&t
\end{pmatrix},$ be $N \times N$. Denote by $\theta^{\beta}(j)$ the
following matrix
\begin{equation*}
\theta^{\beta}(j):=
\left(
\begin{array}{c|ccc}
1_\beta&&&\\
\hline
&0&1_j&\\
&1_j&0&\\
&&&1_\infty
\end{array}
\right)\in\G(\beta).
\end{equation*}
Consider the sequence
$$
\G(\alpha)\cdot g_1\theta^{\beta}(j)\,g_2\cdot\G(\gamma)
\, \in\, \G(\alpha) \setminus \G/\G(\gamma)
.
$$
It is more or less obvious that this sequence is eventually constant and its limit
is
\begin{align}
\label{eq:circ0}
&\frg_1\circ\frg_2=\\
&=
\left[\left(\begin{array}{c|ccc}
a&b&&\\
\hline
c&d&&\\
&&1_{L}&\\
&&&1_\infty
\end{array}\right)
\left(
\begin{array}{c|ccc}
1_\beta&&&\\
\hline
&0&1_L&\\
&1_L&0&\\
&&&1_\infty
\end{array}
\right)
\left(\begin{array}{c|ccc}
p&q&&\\
\hline
r&t&&\\
&&1_{L}&\\
&&&1_\infty
\end{array}\right)\right]_{\alpha\gamma},
\nonumber
\end{align}
where $L\ge N-\beta$. The final expression is
\begin{equation}
\frg_1\circ\frg_2=
\left[
\begin{array}{c|ccc}
ap&aq&b&\\
\hline
cp&cq&d&\\
r&t&0&\\
&&&1_\infty
\end{array}
\right]_{\alpha\gamma}
\sim\,\,
\left[
\begin{array}{c|ccc}
ap&b&aq&\\
\hline
cp&d&cq&\\
r&0&t&\\
&&&1_\infty
\end{array}
\right]_{\alpha\gamma}.
\label{eq:circ}
\end{equation}
{\it In calculations below we use the last expression for $\circ$-product.}

It is is easy to verify that this multiplication is associative, i.e., for any
$$
\frg_1 \in \G(\alpha) \setminus \G/\G(\beta),\quad
\frg_2 \in \G(\beta) \setminus \G/\G(\gamma),\quad
\frg_3 \in \G(\gamma) \setminus \G/\G(\delta),
$$
we have
$$(\frg_1 \circ \frg_2) \circ \frg_3 = \frg_1 \circ (\frg_2 \circ \frg_3).
$$
In other words, we get a  category. Objects of this category are numbers 
$\alpha = 0$,
1, 2, \dots. Sets of morphisms are
$$\Mor(\beta,\alpha):=\G(\alpha) \setminus \G/\G(\beta).
$$
The multiplication is given by formula 
(\ref{eq:circ}). Denote this {\it category} by $\cK$.

\sm 

The {\it group of automorphisms} $\Aut_\cK(\alpha)$
 is $\GL(\alpha,\ZZ)$, it consists
of double cosets of the form 
$\left[\begin{array}{c|c}
a&0\\
\hline
0&1_\infty
\end{array}\right]
$.

\sm

Next, the map $g \mapsto g^{-1}$ induces maps
$$
\G(\alpha) \setminus \G/\G(\beta)\to \G(\beta) \setminus \G/\G(\alpha),
$$
denote these maps by $\frg\mapsto\frg^*$.
It is easy to see that  we get an {\it involution in the
category} $\cK$, i.e.,
$$(\frg_1 \circ\frg_2)^* = \frg_2^* \circ \frg_2^*.$$

\sm

 The map $g \mapsto g^{t-1}$ determines an {\it automorphism of the category} $\cK$, denote it by $\frg\mapsto \frg^\bigstar$. It sends objects to themselves and
 $$
 (\frg_1\circ\frg_2)^\bigstar=\frg_1^\bigstar\circ\frg_2^\bigstar.
 $$
 
 {\sc Remarks on notation.}
 1) In formulas (\ref{eq:g1g2}), (\ref{eq:circ0}), (\ref{eq:circ}),
  the last columns, the last rows, and the blocks
 $1_\infty$ contain no information and only enlarge sizes of matrices. 
 For this reason,
 below we will omit them. 
 Precisely, for a matrix $ \begin{pmatrix}
 a&b\\c&d
 \end{pmatrix}$
 of finite size we denote
 $$
 \left[
 \begin{array}{c|c}
 a&b\\
 \hline c&d_\star
 \end{array}\right]:=
  \left[
 \begin{array}{c|cc}
 a&b&0\\\hline c&d&0
 \\
 0&0&1_\infty
 \end{array}\right]
 \qquad
 \left(
 \begin{array}{c|c}
 a&b\\
 \hline c&d_\star
 \end{array}\right):
 =
 \left(
 \begin{array}{c|cc}
 a&b&0\\
 \hline
 c&d&0
 \\
 0&0&1_\infty
 \end{array}\right)
 .
 $$

2) We will denote a multiplication of $[g]$ by an automorphism 
$A$
as $A\cdot [g] $,
\begin{align*}
A\cdot \left[\begin{array}{c|c}
a&b\\\hline c&d_\star
\end{array}\right]&:=
\left[\begin{array}{c|c}
A&0\\
\hline
0&1
\end{array}
\right]\circ 
\left[\begin{array}{c|c}
a&b\\
\hline
c&d_\star
\end{array}
\right]=\left[\begin{array}{c|c}
Aa&Ab\\
\hline
c&d_\star
\end{array}
\right];
\\
\quad\quad
\left[\begin{array}{c|c}
a&b\\\hline c&d_\star
\end{array}\right]\cdot A'
&
:=
\left[\begin{array}{c|c}
a&b\\\hline c&d_\star
\end{array}\right]
\circ
\left[\begin{array}{c|c}
A'&0\\\hline 0&1
\end{array}\right]
=
\left[\begin{array}{c|c}
aA'&b\\\hline cA'&d_\star
\end{array}\right]. 
\qquad\quad\quad
\boxtimes 
\end{align*}

 \sm 
 
 {\bf\punct The multiplicativity theorem.%
 	\label{ss:multiplicativity}} Consider a unitary representation $\rho$ of
 the group $\G$ in a Hilbert space $H$. Denote by $H_\alpha \subset H $
 the space of $\G(\alpha)$-fixed
 vectors. Denote by $P_\alpha$ the operator of orthogonal projection to $H_\alpha$.
 
 \begin{proposition}
 	\label{pr:theta}
 	{\rm a)} For any $\beta$ the sequence 
 	$\rho\bigl(\theta^{\beta}(j)\bigr)$ converges to $P_\beta$
 	in the weak operator topology.
 	
 	\sm
 	
 	{\rm b)} The space $\cup H_\alpha$ is dense in $H$.
 \end{proposition}

The first statement is Lemma 1.1 from
\cite{Ner-p}, the claim b) is a special case of
Proposition VII.1.3 from \cite{Ner-book}.
\hfill $\square$ 

\sm

Let $g \in \G$, $\alpha$, $\beta \in \Z_+$. Consider the operator
$$
\wt \rho_{\alpha\beta}(g):H_\beta\to H_\alpha
$$
given by
$$
\wt \rho_{\alpha\beta}(g): = P_\alpha \rho(g)\Bigr|_{H_\beta}.
$$
It is easy to see that for $h_1\in \G_\alpha$, $h_2 \in \G_\beta$ we have
$$
\wt \rho_{\alpha\beta}(g)=\wt \rho_{\alpha\beta}(h_1gh_2),
$$
i.e., $\wt \rho_{\alpha\beta}(g)$ actually depends on the double coset $\frg$ containing $g$.

\begin{theorem}
{\rm a)} The map $g \mapsto \wt \rho_{\alpha\beta}(\frg)$ 
is a representation of the category $\cK$, i.e.,
for any $\alpha$, $\beta$, $\gamma$ for any
 $\frg_1 \in \Mor(\beta,\alpha)$, $\frg_2 \in \Mor(\gamma, \beta)$ we have
 $$
 \wt\rho_{\alpha\beta}(\frg_1)\, \wt \rho_{\beta\gamma}(\frg_2)=
 \wt \rho_{\alpha\gamma}(\frg_1\circ\frg_2).
 $$
 
{\rm b)} $\wt \rho$ is a $*$-representation, i.e.,
$$
\wt\rho_{\alpha\beta}(\frg)^*=\wt\rho_{\beta\alpha}(\frg^*).
$$
\end{theorem}

The statement a) is an automatic corollary of Proposition
\ref{pr:theta}, see
\cite{Ner-p}, Theorem
2.1, the statement b) is obvious.
\hfill $\square$

\sm

{\sc Remark.}
 The considerations
of Subsections \ref{ss:cosets}, \ref{ss:multiplicativity} are one-to-one repetitions of similar statements for
real classical groups and symmetric groups, see
\cite{Olsh}, \cite{Olsh-symm},
\cite{Ner-sph},
\cite{Ner-symm}. Further considerations drastically differ from these theories.
\hfill $\boxtimes$

\sm

{\bf \punct Structure of the paper.}
We  derive the classification of unitary representations of
$\G$ from the multiplicativity theorem and the following argumentation. 
The semigroups $\Gamma(m):=\End_\cK(m)$ are finite. 
It is known that a finite semigroup with an involution
has a faithful $*$-representation in a Hilbert space if and only if
it is an inverse semigroup (see discussion below, Subs. \ref{ss:inverse}).
 More generally, if a category having finite sets of morphisms  acts
 faithfully
in Hilbert spaces, then it must be an inverse category, see
\cite{Kas}. However, semigroups
$\End_\cK(\alpha)$ are not inverse%
\footnote{This was observed by Olshanski \cite{Olsh}
	for $\GL(\infty,\F_p)$.}, and $*$-representations of $\cK$ pass through a smaller
category.

Section \ref{s:pr} contains preliminary remarks on
inverse semigroup and construction of an inverse category $\cL$, which
is a quotient of $\cK$. This provides us lower estimate of  maximal
inverse semigroup quotients of semigroups $\Gamma(m)$.

In Section \ref{s:inv} we examine idempotents in  maximal inverse semigroup
quotients $\inv(\Gamma(m))$ of $\Gamma(m)$. In Section \ref{s:red}  we show that some of 
idempotents of $\inv(\Gamma(m))$ act by the same operators in all representations
of $\G$. Next, for any representation
of $\G$ there is a minimal $m$ such that $H_m\ne 0$. 
In Section \ref{s:red-m}
   we examine the image of $\Gamma(m)$ in such representation.
   
In   Section \ref{s:bullet} we discuss  properties of the groups
$\G^\circ[L;M]$ and $\G^\bullet[L;M]$.
    
The final part of the proof is contained in Section \ref{s:fin}.

\sm

\section{The reduced category and inverse semigroups%
\label{s:pr}}

\COUNTERS

{\bf \punct Notation.} Below we work only with the group $\G:=\GL(\infty, \ZZ)$.
To simplify notation, we write
$$
\GL(m):=\GL(m,\ZZ),\qquad \Gamma(m):=\End_\cK(m), \qquad \frl^m:=(\ZZ)^m.
$$

For a unitary representation $\rho$ of a
$\G$ we define the {\it height} $h(\rho)$ as the minimum of $\alpha$ such that $H_\alpha\ne 0$.

By $x(\md p)$ we denote a reduction of an object
(a scalar, a vector, a matrix) defined over $\ZZ$ modulo $p$, i.e. to the field
$\F_p$. Notice that a square matrix $A$ of finite size over $\ZZ$ is invertible if 
and only if $A(\md p)$ is invertible. A matrix $B$ is nilpotent
(i.e., $B^N=0$ for sufficiently large $N$) if and only if $B(\md p)$
is nilpotent.

We use several symbols for equivalences in $\Mor_\cK(\beta,\alpha)$,
 the $\sim$ was defined by (\ref{eq:sim}),
the symbols
$$
\equiv,\qquad\thickapprox,\qquad\thickapprox_m
$$ 
are defined in the next two subsections.

\sm

{\bf\punct The reduced category $\red(\cK)$.} 
Let $\frg_1$, $\frg_2\in \Mor(\beta,\alpha)$.
We say that they are {\it $\thickapprox$-equivalent}
if for any unitary representation
of $\G$ we have $\wt\rho_{\alpha\beta}(\frg_1)=\wt\rho_{\alpha\beta}(\frg_2)$.
The {\it reduced category} $\red(\cK)$ is the category, whose objects
are nonnegative integers and morphisms $\beta\to\alpha$ are
$\thickapprox$-equivalence classes of $\Mor(\beta,\alpha)$. Denote by $\red(\Gamma(m))$
semigroups of endomorphisms of $\red(\cK)$.

Also we define a weaker equivalence, $\frg_1\thickapprox_m \frg_2$
if $\wt \rho_{\alpha\beta}(\frg_1)=\wt \rho_{\alpha\beta}(\frg_1)$
for all $\rho$ of height $\ge m$. Denote by 
$\red_m(\cK)$ the corresponding {\it $m$-reduced category}.

\sm

 Our proof of Theorem \ref{th:main} is based
on an examination of the categories $\red(\cK)$ and $\red_m(\cK)$.
We obtain  an information sufficient for a classification of representations
of $\G$. However, the author does not know an answer to the following question.

\begin{question}
Find a transparent
 description of the category $\red(\cK)$.
	\end{question}

{\bf\punct  Inverse semigroups.%
\label{ss:inverse}}
Let $\cP$ be a {\bf finite} semigroup with an involution $x\mapsto x^*$. Then the following
conditions are equivalent.

\sm

A) $\cP$ admits a faithful representation in a Hilbert space.

\sm

B) $\cP$ admits an embedding to a semigroup of partial bijections%
\footnote{Recall that a {\it partial bijection} $\sigma$
from a set $A$ to a set $B$ is a bijection 
from a subset $S$ of $A$ to a subset $T$ of $B$, see e.g., \cite{Law}
or \cite{Ner-book}, Sect. VIII.1. The adjoint partial bijection
$\sigma^*:B\to A$ is the inverse bijection $T$ to $S$.} of a finite set
compatible with the involutions in $\cP$ and in partial bijections.

\sm

C) $\cP$ is an {\it inverse semigroup} (see \cite{CP}, \cite{Law}, \cite{Kur}),
i.e., for any $x$ we have
\begin{equation}
xx^*x = x, \qquad x^*xx^* = x^*
\label{eq:xxx}
\end{equation}
and any two idempotents in $\cP$ commute.

\sm

Discuss briefly some properties if inverse semigroups.
Any idempotent in $\cP$ is self-adjoint, and for any $x$, the element
$x^*x$ is an idempotent.
Since idempotents commute, a product of idempotents is an idempotent.
The semigroup of idempotents has a natural {\it partial order},
$$
x\preceq y \qquad \text{if} \qquad xy=x.
$$
We have $xy\preceq x$. If $x\preceq y$ and $u\preceq v$,
then $xu\preceq yv$.
Since our semigroup is finite, the product of all idempotents is a {\it minimal idempotent}
$\0$, we have $\0 x=x\0=\0$ for any $x$.

\sm

Let $\cR$ be a finite semigroup with involution. Then there exists an inverse semigroup $\inv(\cR)$ and epimorphism  $\pi:\cR \to \inv(\cR)$ such that any homomorphism
$\psi$ from $\cR$ to an inverse semigroup $\cQ$ has the form $\psi=\kappa\pi$ for some homomorphism
$\kappa: \inv(\cR) \to\cQ$. We say that $\inv(\cR)$ is the {\it maximal inverse semigroup quotient}
of $\cR$.

\sm

\begin{lemma}
The semigroups $\Gamma(m)$ are finite.
\end{lemma}

This is a corollary of the following statement, see  \cite{Ner-p}, Lemma 4.1.a. 

\begin{lemma}
	\label{l:3m}
Any double coset in $\G(m)\setminus \G/\G(m)$ has a representative
in 	$\GL(3m)$.
\end{lemma}

\sm

 We consider the following quotients of $\Gamma(m)$:

\sm 

1) $\inv(\Gamma(m))$ is the maximal inverse semigroup quotient of $\Gamma(m)$;

\sm 

2) $\red(\Gamma(m)):=\End_{\red \cK}(m)$;	

\sm 

3) $\red_m (\Gamma(m)):=\End_{\red_m(\cK)}(m)$.

\sm

We have the following sequence of epimorphisms%
\footnote{All these semigroups are different.}:
$$
\Gamma(m)\to \inv(\Gamma(m))\to \red(\Gamma(m))\to \red_m(\Gamma(m)).
$$
 For $g\in \G_\fin$
we denote by $[g]_{mm}$ the corresponding element of $\Gamma(m)$
and by $[[g]]_{mm}$ the corresponding element of $\inv\bigl(\Gamma(m)\bigr)$.
The equality in $\Gamma(m)$ we denote by $\sim$, in $\inv\Gamma((m))$ by
$\equiv$, in $\red(\Gamma(m))$ by $\thickapprox$, in $\red_m(\Gamma(m))$ by $\thickapprox_m$. Denote by $[[g_1]]\di [[g_2]]$ the product in 
$\inv(\Gamma(m))$.

\sm

Our next purpose is to present some (non-maximal) inverse semigroup quotient
of $\Gamma(m)$.

\sm

{\bf\punct  The category $\cL$ of partial isomorphisms.}
Let $V$, $W$ be modules over $\ZZ$. A {\it partial isomorphism} $p : V \to W$
 is
an isomorphism of a submodule $A \subset V$ to a submodule $B \subset W$. We denote
$\dom p := A$, $\im p := B$. By $p^*$ we denote the inverse map $B\to A$. Let
$p : V \to W$, $q : W \to Y$ be partial isomorphisms. Then the {\it product} $pq$ is defined
in the following way:
$$\dom pq := p^*(\dom q)\cap \dom p,$$
for $v \in \dom pq$ we define $v(pq) = (vp)q$.

A partial isomorphism $p$ is an {\it idempotent} if $\dom p = \im p$
 and $p$
is an identical map.

\sm

Objects of the category $\cL$ are modules
$$\frl_+^\alpha \oplus \frl_-^\alpha:=
(\ZZ)^\alpha\oplus (\ZZ)^\alpha$$
equipped with the following pairing
$$\{v_+; v_-\} := \sum_j v_+^j v_-^j=v_+(v_-)^t,$$
 where $v_\pm \in  \frl_\pm^\alpha$.
We say that two partial isomorphisms
$$\xi_+ : \frl_+^\alpha \to \frl_+^\beta,
\qquad
\xi_- : \frl_-^\alpha \to \frl_-^\beta$$
are {\it compatible} if for any $y_+\in\dom \xi_+$ and
 $y_-\in\dom \xi_-$, we have
$$
\{\xi_+(y_+),\xi_-(y_-) \}=\{y_+,y_-\}.
$$
Next, we define a {\it category} $\cL$. Its objects are spaces 
$\frl_+^\alpha\oplus \frl_-^\alpha$ and
morphisms are pairs of compatible partial isomorphisms
$\xi_+ : \frl_+^\alpha \to \frl_+^\beta$,
$\xi_- : \frl_-^\alpha \to \frl_-^\beta$.

The category $\cL$ is equipped with an involution
$$(\xi_+,\xi_-)^*=(\xi_+^*,\xi_-^*)$$
and an automorphism
$$
(\xi_+,\xi_-)^\bigstar=(\xi_-,\xi_+).
$$

\begin{lemma}
	The semigroups $\End_\cL(m)$ are inverse.
\end{lemma}

Indeed, $\End_\cL(m)$ is a semigroup of partial bijections 
of a finite set $\frl^m_+\oplus \frl^m_-$. The whole category
$\cL$ is inverse for the same reason.
\hfill $\square$

\sm

{\bf\punct The functor $\Pi : \cK \to\cL$.%
\label{ss:KL}} Consider $g \in \G_\fin$.
 Let actually $g$ be
contained in $\GL(N)$. Represent $ g$ as a block
$(\beta+(N-\beta))\times (\alpha+(N-\alpha))$
matrix and $g^{-1}$ as an $(\alpha+(N-\alpha))\times(\beta+(N-\beta))$-matrix,
$$
g=\begin{pmatrix}
a&b\\c&d_\star
\end{pmatrix}, \qquad
g^{-1}=\begin{pmatrix}
A&B\\C&D_\star
\end{pmatrix}.
$$
Define maps
$\xi_\pm:\frl_\pm^\alpha\to\frl_\pm^\beta$
  by:
  
  \sm
  
$\bullet$ $\dom \xi_+ := \ker b$ and $\xi_+$ is the restriction of $a$ to $\ker b$;

 \sm

$\bullet$  $\dom \xi_- := \ker C^t$ and $\xi_-$ is the restriction of $A^t$ to $\ker C^t$.

\begin{proposition}
{\rm a)} The pair $\xi_+$, $\xi_-$ depends only on the double coset containing $\frg$.

 \sm
 
 {\rm b)} Partial isomorphisms $\xi_+$, $\xi_-$ are compatible.

\sm

{\rm c)} The map $\frg\mapsto(\xi_+, \xi_-)$ determines a functor from the category $\cK$ to the
category $\cL$.
\end{proposition}	
	
Denote this functor by $\Pi$. By $\Pi(\frg)$ we denote the morphism of $\cL$ corresponding to $\frg$.
 We have
\begin{equation}
\Pi(\frg^*)=(\Pi(\frg))^*,\qquad
\Pi(\frg^\bigstar)=(\Pi(\frg))^\bigstar.
\label{eq:*-star}
\end{equation}

{\sc Proof.} 
For any invertible matrix $v$ we have, $\ker b = \ker bv$.
 Therefore $\xi_+$  depends
only on a double coset. For $\xi_-$ we apply (\ref{eq:*-star}).

\sm

b) Let $v \in \ker b$, $w \in \ker C^t$. Then
$$
\{v,w\}=vw^t= v(aA + bC)w^t = va · (wA^t)^t + vb · (wC^t)^t = \{va, wA^t\}+0
.
$$

c) We look to formula (\ref{eq:circ})
 for a product in $\cK$. The new $\xi_+$ is a restriction
of $ap$ to $\ker b \cap \ker aq$. This is the product of two $\xi$-es.
\hfill $\square$

 \sm
 
 {\sc Remark.}
 According Olshanski \cite{Olsh}, for the case
 $\GL(\infty,\F_p)$ the functor $\Pi:\cK\to\cL$ determines
 an isomorphism of categories $\red(\cK)\to \cL$.
 However, for
 $\mu>1$ the maps $\Pi: \red(\Gamma(m))\to \Mor_\cL(m)$ are neither surjective, nor injective.
 However we will observe, that $\Pi$ induce  isomorphisms  of  semigroups of idempotents; this provides us an important argument
for the proof of Proposition \ref{pr:structure}.
 \hfill $\boxtimes$


	
	
%
%

\section{Idempotents in $\inv(\Gamma(m))$%
\label{s:inv}}

\COUNTERS

Here we examine idempotents in the semigroup $\inv(\Gamma(m))$. The main statement
of the section is Proposition \ref{pr:idempotents-4}.

\sm 

{\bf\punct Projectors%
\footnote{This subsection contains generalities, $\cK$ is an ordered category
in the sense of \cite{Ner-book}, Sect. III.4, this implies all statements of the subsection.}
$P_\alpha$.}
Consider an irreducible representation $\rho$ of $\G$, let subspaces $H_m\subset H$ and orthogonal projectors $P_m:H\to H_m$ be as above.

\begin{lemma}
	\label{l:P-alpha}
{\rm a)}	The projector 
$$P_\alpha\Bigr|_{H_m}:H_m\to H_\alpha$$
is given 
	by the operator
	$\wt\rho_{mm}\bigl(\Theta_{[m]}^{\alpha}\bigr)$, where 
	\begin{equation}
	\Theta_{[m]}^\alpha:=
	\left[
	\begin{array}{cc|cc}
	1_\alpha&0&0&0\\
	0&0&1_{m-\alpha}&0\\
	\hline
	0&1_{m-\alpha}&0&0\\
	0&0&0&1_\infty
	\end{array}
	\right]_{mm}\in\Gamma(m).
	\end{equation}
	
	{\rm b)} The tautological embedding $H_\alpha\to H_m$
	is defined by the operator $\wt\rho_{m\alpha}\bigl(\Lambda^\alpha_{[m]}\bigr)$, 
	there
	\begin{equation*}
	\Lambda^\alpha_{[m]}:=
	\left[\begin{array}{cc|cc}
	1_\alpha&0&0\\
	\hline 
	0&1_{m\alpha}&0\\
	0&0&1_\infty
	\end{array}\,\,
	\right]_{\alpha m}
	\in \Mor_\cK(\alpha,m).
	\end{equation*}
	
	{\rm c)} The orthogonal projector $H_m\to H_\alpha$
	is given by $\wt\rho_{\alpha m}\bigl((\Lambda_{[m]}^\alpha)^*\bigr)$
	$$
	(\Lambda^\alpha_m)^*:=
	\left[\begin{array}{c|cc}
	1_\alpha&0&0\\
	0&1_{m-\alpha}&0\\
	\hline 
	0&0&1_\infty
	\end{array}\right]_{m\alpha}
	\in \Mor_\cK(m,\alpha).
	$$
\end{lemma}

{\sc Proof.} a) We apply Proposition
\ref{pr:theta}.a.
 For $j>m-\alpha$  we have $[\theta^{\alpha}(j)]_{mm} =\Theta_{[m]}^\alpha$.
 The same argument proves b) and c).
\hfill $\square$

\begin{lemma}
	\label{l:embed}
{\rm a)} The map
	$$
	\iota^\alpha_m:
	\left[\begin{array}{c|c}
		a&b\\\hline c&d_\star
	\end{array}\right]_{\alpha\alpha}
	\mapsto\left[\begin{array}{cc|cc}
		a&0&b&0\\
		0&0&0&1_{m-\alpha}
		\\\hline c&0&d&0\\
		0&1_{m-\alpha}&0&0_\star
	\end{array}\right]_{mm}$$
	is a homomorphism $\Gamma(\alpha)\to \Gamma(m)$.
	
	\sm 
	
	{\rm b)} We have
	$$
	\iota^\alpha_m(\frg)\sim\Lambda^\alpha_n \circ \frg\circ \left(\Lambda^\alpha_n \right)^*.
	$$
\end{lemma}

This follows  from  a straightforward calculation.
\hfill $\square$

\begin{corollary}
	\label{cor:iota}
 The map $\iota_m^\alpha$ is compatible
with representations $\wt \rho$ of $\Gamma(\alpha)$ and $\Gamma(m)$.
Namely, operators $\wt\rho_{mm}(\iota^\alpha_m(\frg))$ have the following
block structure with respect to the decomposition
$H_m=H_\alpha\oplus (H_m\ominus H_\alpha)$:
$$
\wt\rho_{mm}(\iota^\alpha_m(\frg))=
\begin{pmatrix}
\wt\rho_{\alpha \alpha}(\frg)&0\\0&0
\end{pmatrix}.
$$
\end{corollary}

\sm

{\bf\punct Idempotents in $\inv(\Gamma(m))$.}
Here we formulate several lemmas (their proofs occupy  Subsections \ref{ss:4.3}--\ref{ss:4-last}),
as a corollary we get Proposition \ref{pr:idempotents-4}.

\begin{lemma}
	\label{l:Gamma-circ}
Let for	$$[g]=\begin{bmatrix}
		a&b\\c&d_\star
	\end{bmatrix}_{mm}\in \Gamma(m)
	$$
one of the blocks $a$, $d$ be degenerate. Then 
$[[g]]_{mm}\in \inv(\Gamma(m))$ has a representative 
$[g']$, for which both blocks $a$, $d$ are degenerate.	
\end{lemma}

Denote by 
$$\Gamma^\circ(m)$$
the subsemigroup
in $\Gamma(m)$ consisting of all $[g]$, for which both blocks $a$, $d$ are nondegenerate.

\begin{lemma}
	\label{l:idempotents-1}
Any idempotent in $\inv (\Gamma(m))$ has a representative of the form $q\cdot[[R]]\cdot q^{-1}$
with $q$ ranging $\GL(m)$ and  $R$ having the form
\begin{equation}
[R]:=
\left[
\begin{array}{cc|cc}
1_\alpha&0&\phi&0\\
0&0&0&1_{m-\alpha}\\
\hline
\psi&0&\kappa&0\\
0&1_{m-\alpha}&0&0_\star
\end{array}
\right]_{mm}\in\Gamma(m),
\label{eq:[R]}
\end{equation}
 where
$$
\left[
\begin{array}{c|c}
1&\phi\\
\hline
\psi&\kappa_\star
\end{array}
\right]_{\alpha\alpha}\in \Gamma(\alpha)
$$
represents an idempotent in $\inv(\Gamma^\circ(\alpha))$. The parameter $\alpha$ ranges in the set 
$0$, $1$, $2$, \dots, $m$.
\end{lemma}	

{\sc Remark.}
Denote
$$
R^\square:=
\left[
\begin{array}{cc|c}
1_\alpha&0&\phi\\
0&1_{m-\alpha}&0\\
\hline
\psi&0&\kappa_\star
\end{array}
\right]_{mm}.
$$
Then the following elements of $\Gamma(m)$ coincide:
\begin{equation}
R=R^\square\Theta_m^\alpha=\Theta_m^\alpha R^\square =\Theta_m^\alpha R^\square \Theta_m^\alpha.
\label{eq:RthetaR}
\end{equation}

Denote
$$
X(b,c):=\left(
\begin{array}{c|cc}
1_m&b&0\\
\hline
0&1&0\\
c&0&1_\star
\end{array}
\right)\in\G_\fin.
$$

\begin{lemma}
	\label{l:idempotents-2}
Elements of the form $[X(b,c)]$	
are idempotents in $\Gamma^\circ(m)$.
They depend only on $\ker b$ and $\ker c^t\subset \frl^m$.
\end{lemma}	

Let $L:=\ker b$ and $M:=\ker c^t$. Denote
\begin{equation}
\cX[L,M]:=[X(b,c)].
\label{eq:cX}
\end{equation}

\begin{lemma}
\label{l:idempotents-21}
We have 
$$
\cX[L_1,M_1]\,\,\cX[L_2,M_2]=\cX[L_1\cap L_2, M_1\cap M_2].
$$	
\end{lemma}

\begin{lemma}
	\label{l:idempotents-3}
	Any idempotent in $\inv(\Gamma^\circ(m))$ has the form 
	$\cX[L,M]$.
\end{lemma}

\begin{corollary}
	\label{cor:} Idempotents $\cX[L,M]$ are pairwise distinct in $\inv(\Gamma^\circ(m))$.
\end{corollary}

{\sc Proof.} Indeed, $\End_\cL(m)$ is an inverse semigroup, therefore
we have a chain of maps
$$
\Gamma^\circ(m)\to
\inv(\Gamma^\circ(m))\to \inv(\Gamma(m))\to \Mor_\cL(m).
$$ 
The image of $X(b,c)$ in $\Mor_\cL(m)$ is precisely the pair of identical
partial isomorphisms $M\to M$, $L\to L$.  Therefore for nonequivalent $X(b,c)$
we have different images.
\hfill $\square$

\begin{proposition}
	\label{pr:idempotents-4}
	Any idempotent in $\inv(\Gamma(m))$ has a representative of the form
	$$
	q\cdot
	\left[\begin{array}{cc|ccc}
	1_{\alpha}&0&b&0&0\\
	0&0&0&0&1_{m-\alpha}\\
	\hline
	0&0&1&0&0\\
	c&0&0&1&0\\
	0&1_{m-\alpha}&0&0&0_\star
	\end{array}\right]_{mm}\cdot
	q^{-1},
	$$
	where $q\in \GL(m)=\Aut_\cK(m)$.
\end{proposition}

{\sc Proof.} Lemma \ref{l:embed} defines a canonical embedding
$i^\beta_m:\Gamma(\beta)\to \Gamma(m)$ for $\alpha<m$. 
By Lemma \ref{l:idempotents-1}
any idempotent in $\inv(\Gamma(m))$ is equivalent to
and idempotent lying in some $i^\alpha_m(\Gamma^\circ(\alpha))$.
Lemma \ref{l:idempotents-3} gives us a canonical form of this idempotent.
\hfill $\square$

\sm 

Now we start proofs of Lemmas \ref{l:Gamma-circ}--\ref{l:idempotents-3},

\sm

 {\bf\punct Proof of Lemma \ref{l:Gamma-circ}.%
 \label{ss:4.3}}
Clearly $\Gamma(m)\setminus \Gamma^\circ(m)$
is a two-side ideal in $\Gamma(m)$.
Since $[[g\circ(g^{-1}\circ g)]]_{mm}=[[g]]_{mm}$,
it is sufficient to prove the statement for  idempotents.

Let 
\begin{equation}
\label{eq:gg}
g=\begin{pmatrix} a&b\\c&d_\star\end{pmatrix},\qquad g^{-1}=:\begin{pmatrix} A&B\\C&D_\star\end{pmatrix}.
\end{equation}
Then
$$
[[g]]\di [[g]]^*\equiv [[g\circ g^{-1}]]\equiv
\left[\left[ \begin{array}{c|c} aA&*\\
\hline
*&*
\end{array}\right]\right].
$$
If $a$ is degenerate, then $aA$ is degenerate. Now let $a$ be non-degenerate, $d$ degenerate. 
Since the matrices (\ref{eq:gg}) are inverse one to another, we have
$$
aA=1-bC, \qquad Dd=1-Cb.
$$
We see that $(1-Cb) (\md p)$ is degenerate, $(1-bC)(\md p)$ also is degenerate, and therefore $aA$ is degenerate.

\sm

{\bf\punct Proof of Lemma \ref{l:idempotents-1}.}
{\sc Step 1.} 

\begin{lemma}
	\label{l:[[u]]}
	Let $x$ be an idempotent in $\inv(\Gamma(m))$. Then it can be represented as
	$[[u]]$, where $u=u^{-1}$.
\end{lemma}

{\sc Proof.}  
Let $x=[[g]]$. Then
$$
x=
[[g]]\di [[g]]^*= [[g\circ g^{-1}]]=
[[g \theta^{m}(j) g^{-1}]]
$$
for sufficiently large $j$. We set $u:=g \theta^{m}(j) g^{-1}$.
\hfill $\square$

\begin{lemma}
	\label{l:ggg}
	Let $g=g^{-1}\in \G_\fin$. For any $N>0$ there exists a representative $r\in \G_\fin$ of $[g]^{\circ 2N}$
	such that $r=r^{-1}$.
\end{lemma}

{\sc Proof.} Let actually $g\in\GL(m+l)$. Then we choose the following
 representative of $[g]^{\circ 8}$:
$$
r=g\,\,\theta^{m}(l)\,\, g\,\, \theta^{m}(2l)\,\, g\,\, \theta^{m}(4l)\,\, g \,\,\theta^{m}(8l)
\,\,g \theta^{m}(4l)\,\,g \,\,\theta^{m}(2l)\,\, g\,\, \theta^{m}(l)\, \,g.
$$

\sm

{\sc Step 2.} 

\begin{lemma}
	\label{l:gN}
	Let $g=g^{-1}=\begin{pmatrix} a&b\\c&d_\star\end{pmatrix}\in \G$.
	Then there exists a matrix
	$$Z=\left(
	\begin{array}{c|c}
	\zeta&0\\\hline
	0&1_\star
	\end{array}
	\right)\in \Aut_\cK(m), \qquad \text{where  $\zeta\in \GL(m)$,}
	$$
	and $N$ such that
	$$
	[[Z\cdot g\cdot Z^{-1} ]]^{\di N}=
	\left[\left[ \begin{pmatrix} \zeta &0\\0&1_\star\end{pmatrix}
	\begin{pmatrix} a&b\\c&d_\star\end{pmatrix}
	\begin{pmatrix} \zeta &0\\0&1_\star\end{pmatrix}^{-1}\right]\right]^{\di N}
	$$
	has a form
	$$ r=
	\left[\left[  \begin{array}{cc|c}
	0&0&*\\
	0&1_k&*\\
	\hline
	*&*&*
	\end{array}\right]\right],
	$$
	where $k$ is the rank of the reduced matrix $a^m (\md p)$.
\end{lemma}

Clearly our lemma is a corollary of the following statement:

\begin{lemma}
	For any $m\times m$ matrix $a$
	over $\ZZ$ there exists  $\zeta\in \GL(m)$ and $N$
	such that 
	$$
	(\zeta a \zeta^{-1})^N=\begin{pmatrix}
	0&0\\0&1_k
	\end{pmatrix}.
	$$
\end{lemma}

{\sc Proof.} 
We split the operator $a(\md p)$ over the field $\F_p$ as a direct sum of a nilpotent 
part $S$ and a invertible part $T$.
For sufficiently large $M$ the matrix
$\begin{pmatrix}
S&0\\0&T
\end{pmatrix}^M$ has 
the form $\begin{pmatrix}
0&0\\0&P
\end{pmatrix}$ with a nondegenerate $P$. Since the group $\GL(k,\F_p)$
is finite,  $P^L=1_k$ for some $L$.

Thus without loss of generality, we can assume that $a$ has a form
$$
a=\begin{pmatrix}
p \alpha &p\beta\\
p\gamma& 1+p\delta
\end{pmatrix},
$$
where $\alpha$, $\beta$, $\gamma$, $\delta$ are matrices over $\ZZ$. We conjugate it as follows
$$
\begin{pmatrix}
1&pu\\0&1
\end{pmatrix}
\begin{pmatrix}
p \alpha &p\beta\\
p\gamma& 1+p\delta
\end{pmatrix}
\begin{pmatrix}
1&-pu\\0&1
\end{pmatrix}=
\begin{pmatrix}
*& \boxed{-p^2(\alpha u+u\gamma u)+p\bigl(\beta+u (1+p\delta)\bigr)}\,\\
*&*
\end{pmatrix}.
$$
We wish to choose $u$ to make zero in the boxed block. It is sufficient to find a matrix  $u$
satisfying the following equation:
\begin{equation}
u=\bigl(-\beta+p(\alpha u+u\gamma u)\bigr)(1+p\delta)^{-1}
=
-\beta+p(-\delta+\alpha u+u\gamma u)(1+p\delta)^{-1}
.
\label{eq:uu}
\end{equation}
We look for a solution
in the form
$$
u=\sum\nolimits_{k=0}^\mu p^k S_k.
$$
First, we consider $S_k$ as formal noncommutative variables. Then we get a system
of equations of the form
$$
S_0=-\beta,\qquad S_k=F_k(\alpha,\beta, \gamma,\delta; \,S_0,S_1,\dots, S_{k-1}),
$$
where $F_k$ are polynomial expressions with integer coefficients.
These equations can be regarded as	recurrence  formulas for $S_k$.
In this way we get a solution $u$.


\sm 

Thus without a lose of generality we can assume that $a$ has the form 
$$a=\begin{pmatrix}p\alpha'&0\\ p\gamma'& 1+p\delta' \end{pmatrix}.$$
Raising it to $\mu$-th power, we come to a matrix of the form
$$a=\begin{pmatrix}0&0\\ p\gamma''& 1+p\delta'' \end{pmatrix}.$$
We conjugate it as 
$$
\begin{pmatrix}
1&0\\pv& 1
\end{pmatrix}
\begin{pmatrix}0&0\\ p\gamma''& 1+p\delta'' \end{pmatrix}
\begin{pmatrix}
1&0\\-pv& 1
\end{pmatrix}=
\begin{pmatrix}
0&0\\p(\gamma''-(1+p\delta'') v)&(1+\delta'' v)
\end{pmatrix}.
$$
Taking $v=  (1+p\delta'')^{-1} \gamma''$ we kill the left lower block and come to
a matrix of the form
$\begin{pmatrix}0&0\\ 0& 1+p\delta''' \end{pmatrix}
$. Raising it in $p^{\mu-1}$-th power we come to
$\begin{pmatrix}0&0\\ 0& 1 \end{pmatrix}
$. \hfill $\square$

\sm

{\sc Step 3.} 
Thus the element $[[g]]^{\di 2N}$ from Lemma \ref{l:gN}
has a representative of the following
block $(m-k)+k+(m-k)+\infty$ form:
$$r=r^{-1}=
\left(
\begin{array}{cc|cc}
0&0& \beta_{11}&\beta_{12}\\
0&1_k& \beta_{21}&\beta_{22}\\
\hline
\gamma_{11}&\gamma_{12}&\delta_{11}&\delta_{12}\\
\gamma_{21}&\gamma_{22}&\delta_{21}&\delta_{22\star}
\end{array}\right).
$$

\begin{lemma}
	There is a matrix $U=\left(\begin{array}{c|c}
	1_m&0\\
	\hline
	0&u_\star
	\end{array}\right)
	$
	such that $UrU^{-1}$
	has the form
	\begin{equation}
	\wt r=
	\left( 
	\begin{array}{cc|cc}
	0&0&1_{m-k}&0\\
	0&1_k&0&\phi\\
	\hline
	1_{m-k}&0&0&0\\
	0&\psi&0&\kappa_\star
	\end{array}
	\right).
	\label{eq:wtr}
	\end{equation}
\end{lemma}

Recall that $[r]\sim[UrU^{-1}]$.

\sm 

{\sc Proof.}	
Since the matrix $\begin{pmatrix}
\beta_{11}&\beta_{12}
\end{pmatrix}$
is nondegenerate (otherwise $r$ is degenerate), we can choose a conjugation of $r$
by matrices  $U=\begin{pmatrix}
              1_m&0\\0&*
             \end{pmatrix}$
 reducing this block to the form
$\begin{pmatrix}
1&0
\end{pmatrix}$. We have $r^2=1$, evaluating $r^2$  we get $\gamma_{11}$
in the left upper block. Therefore $\gamma_{11}=1$.
Thus we come to new $r$,
$$
r^\sim=
\left(
\begin{array}{cc|cc}
0&0& 1&0\\
0&1_k& \beta_{21}&\beta_{22}\\
\hline
1&\gamma_{12}&\delta_{11}&\delta_{12}\\
\gamma_{21}&\gamma_{22}&\delta_{21}&\delta_{22\star}
\end{array}\right)
$$
with new $\beta$, $\gamma$, $\delta$. Next, we conjugate this matrix by
$$\left(
\begin{array}{c|cc}
1_m&0&0\\
\hline
0&1_{m-k}&0\\
0&-\gamma_{21}&1_\star
\end{array}
\right)
$$
and kill $\gamma_{21}$.
Thus we come to new $r$,
$$
r^{\sim\sim}=
\left(
\begin{array}{cc|cc}
0&0& 1&0\\
0&1_k& \beta_{21}&\beta_{22}\\
\hline
1&\gamma_{12}&\delta_{11}&\delta_{12}\\
0&\gamma_{22}&\delta_{21}&\delta_{22\star}
\end{array}\right)
$$
But $(r^{\sim\sim})^2=1$. Looking to third row and third column of
$(r^{\sim\sim})^2$ we observe that
$$
\text{$\beta_{21}$, $\delta_{11}$, $\delta_{21}$, $\gamma_{12}$, $\delta_{12}$ are zero.}
$$
Thus, $r^{\sim\sim}$
has the desired form.
\hfill $\square$

\sm

{\bf\punct Proof of Lemma \ref{l:idempotents-2}.}
 Denote
$$
[X_+(A)]:=\left[\begin{array}{cc}
1&A\\0&1_\star
\end{array}\right].
$$
We can conjugate this matrix by $\begin{pmatrix}
1&0\\0&u_\star
\end{pmatrix}$.
Therefore a matrix $A$ is defined up to multiplications $A\sim Au$, where $u$
is an invertible matrix.
The invariant of this action is $\ker A$ (this is more or less clear,
formally we can refer to Lemma \ref{l:52} proved below).

Next,
$$
[X_+(A)]\circ [X_+(A)]=\left[\begin{array}{c|cc}
1&A&A\\
\hline 
0&1&0\\
0&0&1_\star
\end{array}\right].
$$
We have $\ker\begin{pmatrix}
A&A
\end{pmatrix}=\ker A$ and therefore $[X_+(A)]$ is an idempotent.
In the same way, $[X_-(B)]:=\left[\begin{array}{cc}
1&0\\B&1_\star
\end{array}\right]$ is an idempotent.
It remains to notice that
$$
[X(A,B)]=[X_+(A)]\di [ X_-(B)].
$$
{\it Thus $[X(A,B)]$ is an idempotent}.

\sm

{\bf\punct Proof of Lemma \ref{l:idempotents-2}.}
In notation of the previous subsection
$$
X_-(A_1)\circ X_-(A_2)\sim X\Bigl( \begin{pmatrix}
A_1&A_2
\end{pmatrix} \Bigr),
$$
i.e., 
$$
\cX[\ker A_1, 0]\di \cX[\ker A_2, 0]\equiv \cX[\ker A_1\cap \ker A_2, 0],
$$
or 
$$\cX[L_1,0]\di\cX[L_2,0]\equiv \cX[L_1\cap L_2,0]$$
On the other hand, we have
$$
\cX[L,0]\di\cX[0,M]\equiv \cX[L,M],
$$
and now the statement becomes obvious.

\sm

{\bf\punct Proof of Lemma \ref{l:idempotents-21}.}
Indeed,
$$
[X(b_1,c_1)]\circ [X(b_1,c_1)]\sim
\left[
X\left(
\begin{pmatrix}
 b_1&b_2
\end{pmatrix},
\begin{pmatrix}
 c_1\\c_2
\end{pmatrix}
\right)
\right]
$$
and $\ker\begin{pmatrix}
 b_1&b_2
\end{pmatrix}=\ker b_1\cap \ker b_2$.
\hfill $\square$

\sm 

{\bf\punct Proof of Lemma \ref{l:idempotents-3}.%
\label{ss:4-last}}

\sm 

{\sc Step 1.}  {\it Any  idempotent $[[g]]\in \inv(\Gamma^\circ(m))$
	 has a representative of
	the form  $\begin{pmatrix}
	1&a\\b&1_\star
	\end{pmatrix}$, where $ab=0$.}

\sm

Let $[[g]]=\left[\left[\begin{array}{cc} \alpha&\beta\\\gamma&\delta_\star\end{array}\right]\right]$
be an idempotent, let $\alpha$, $\delta$ be nondegenerate. By Lemma \ref{l:[[u]]} without loss of generality we can assume $g=g^{-1}$.
Taking an appropriate power $[r]=[g]^{\circ 2N}$, we can achieve $\alpha=1$. By Lemma \ref{l:ggg}, we can assume
$r=r^{-1}$.

Set $r=\begin{pmatrix}1&-a\\b&c_\star\end{pmatrix}$.
Evaluating $r^2=1$ we get the following collection of conditions
$$
\boxed{ab=0}\,,\qquad ac=-a, \qquad cb=-b,\qquad c^2-ba=1.
$$
We replace $r$ by an equivalent matrix
$$r\sim
\begin{pmatrix}1&-a\\b&c_\star\end{pmatrix}\begin{pmatrix}
1&0\\0&c^{-1}_\star \end{pmatrix}=
\begin{pmatrix}1&-ac^{-1}\\b&1_\star\end{pmatrix}   =   \begin{pmatrix}1&a\\b&1_\star\end{pmatrix}  
,                        
$$
here we used the identity $-ac^{-1}=a$.

\sm 

{\sc Step 2.}
We evaluate $[r]^{\circ 2}$,
\begin{multline*}
[r]^{\circ 2}=\left[\left(
\begin{array}{c|cc}1&a&0\\\hline b&1&0\\0&0&1_\star\end{array}\right)
\left( \begin{array}{c|cc}1&0&a\\\hline 0&1&0\\ b&0&1_\star\end{array} \right)\right]=
\left[\begin{array}{c|cc}1&a&a\\\hline b&1&ba\\b&0&1_\star\end{array}\right]
\sim\\\sim\left[\left(\begin{array}{c|cc}1&a&a\\\hline b&1&ba\\b&0&1_\star\end{array}\right)
\left(\begin{array}{c|cc}1&0&0\\\hline 0&1&-ba\\0&0&1_\star\end{array}\right)
\right]=\left[\begin{array}{c|cc}1&a&a-aba\\\hline b&1&0\\b&0&1_\star\end{array}\right].
\end{multline*}
But $ab=0$ and therefore $aba=0$. Repeating the same reasoning, we get
\begin{equation}
[[r]]^{\di N}\equiv[[q]]\equiv
\left[\left[
\begin{array}{c|ccc}
1_m&a&\dots&a\\
\hline
b&1&\dots&0\\
\vdots&\vdots&\ddots&\vdots\\
b&0&\dots&1_\star
\end{array}\right]\right].
\label{eq:aaa}
\end{equation}

{\sc Step 3.}
Next, we set $N=p^\mu$ in formula (\ref{eq:aaa}).
Consider the following block matrix $u$ of size $p^\mu$,
$$
u:=\begin{pmatrix}
1&&&&\\
-1&1&&&\\
0&-1&1&&
\\
\vdots&\vdots&\vdots&\ddots\\
0&0&0&\dots&1
\end{pmatrix},\qquad
u^{-1}:=
\begin{pmatrix}
1&&&&\\
1&1&&&\\
1&1&1&&
\\
\vdots&\vdots&\vdots&\ddots\\
1&1&1&\dots&1
\end{pmatrix}.
$$
We conjugate the matrix $q$ defined by (\ref{eq:aaa}) as
$$
\begin{pmatrix}
1&0\\0&u_\star
\end{pmatrix}
\,q\, 
\begin{pmatrix}
1&0\\0&u^{-1}_\star
\end{pmatrix}.
$$
We have
$$
u\begin{pmatrix}
b\\b\\b\\\vdots
\end{pmatrix}
=\begin{pmatrix}
b\\0\\0\\\vdots
\end{pmatrix},\qquad
\begin{pmatrix}
a&a&a&a&\dots
\end{pmatrix}\, u^{-1}=
\begin{pmatrix}
0&-a&-2a&-3a&\dots
\end{pmatrix},
$$
and we get a matrix of the form $X(A,B)$.



\section{Idempotents in $\red(\Gamma(m))$%
\label{s:red}}

\COUNTERS

Here the main statement is Proposition \ref{pr:coincidence},
which shows that all idempotents in $\red(\Gamma(m))$ have representatives 
in $\red(\Gamma^\circ(m))$, therefore they have the form
$\cX[L,M]$), where $L\subset \frl^m$, $M\subset \frl^m$. The second fact (Proposition \ref{pr:coherence}), 
which is important for the  proof below, is a coherence
of elements $\cX[L,M]$ in different semigroups $\red(\Gamma(n))$.

\sm 

{\bf\punct Coincidence of idempotents.}

\begin{proposition}
	\label{pr:coincidence}
	The following idempotents in $\inv(\Gamma(m))$ coincide as elements of
	$\red(\Gamma(m))$:
\begin{multline*}
[[X_\alpha^\bigcirc(b,c)]]:=
\left[\left[ X\left(\begin{pmatrix}b&0\\0&1_{m-\alpha} \end{pmatrix},
\begin{pmatrix}c&0\\0&1_{m-\alpha} \end{pmatrix}\right)\right]\right]=\\
=
	\left[\left[\begin{array}{cc|cccc}
1_{\alpha}&0&b&0&0&0\\
0&1_{m-\alpha}&0&1_{m-\alpha}&0&0\\
\hline
0&0&1&0&0&0\\
0&0&0&1&0&0\\
c&0&0&0&1&0\\
0&1_{m-\alpha}&0&0&0&1_\star
\end{array}\right]\right]_{mm}
,
\end{multline*}
and
\begin{equation*}
[[X_\alpha^\square(b,c)]]:=
\left[	\left[\begin{array}{cc|ccc}
1_{\alpha}&0&b&0&0\\
0&0&0&0&1_{m-\alpha}\\
\hline
0&0&1&0&0\\
c&0&0&1&0\\
0&1_{m-\alpha}&0&0&0_\star
\end{array}\right]\right]_{mm}.
\end{equation*}	
	\end{proposition}
	
\begin{corollary}
	\label{cor:all-idempotents}
Any idempotent in $\red(\Gamma(m))$ has the form $\cX[L,M]$.
\end{corollary}

{\sc Proof of corollary.}
The semigroup $\red\Gamma(m))$ is a quotient of $\inv(\Gamma(m))$,
the semigroup of idempotents also is a quotient of the semigroup of idempotents.
By Proposition \ref{pr:idempotents-4} all idempotents in $\inv(\Gamma(m))$ have $[[X_\alpha^\bigcirc[b,c]]$.
By Proposition \ref{pr:coincidence}, they also can be written as 
$[[X_\alpha^\square[b,c]]$ .
\hfill $\square$

\sm

Proposition will be proved  below in Subsection \ref{ss:skleika}.

\sm

{\sc Remarks.} a) The idempotents $[[X_\alpha^\bigcirc(b,c)]]$ and $[[X_\alpha^\square(b,c)]]$ are different in $\inv(\Gamma(m))$.
Indeed, we have the following homomorphism from $\Gamma(m)$
to the inverse semigroup $\End_\cL(m)$. On $\Gamma^\circ(m)$
we define it as the map $\Pi$ described in Subs. \ref{ss:KL}.
On the other hand, we send $\Gamma(m)\setminus \Gamma^\circ(m)$ to $\0$,
i.e., to a pair of partial bijections with empty domains of definiteness.
This map separates our idempotents. 

\sm 

b) Idempotents $\cX[L,M]$ are pairwise different in $\red(\Gamma(m))$.
To verify this, consider the representation of $\G$ in $\ell^2(\G/\G[L;M])$.
It is easy to show that $\cX(L,M)$ is the minimal idempotent of
 $\red(\Gamma(m))$ acting in this representation nontrivially.
 \hfill $\boxtimes$

 \sm

{\bf \punct Coherence.}
Let $L$, $M\subset\frl^m$ be submodules. Formula
(\ref{eq:cX}) defines the idempotent $\cX[L,M]=X(b,c)$  as an element
of $\Gamma(m)$, recall that $L=\ker b$, $M=\ker c^t$. However, 
for $n>m$ we can regard $L$, $M\subset \frl^m$ as submodules $L$ in $\frl^n\supset \frl^m$.
In the larger space we have
$$
L=\ker\begin{pmatrix}
b&0\\0&1_{n-m}
\end{pmatrix}, \qquad
M=\ker\begin{pmatrix}
c&0\\0&1_{n-m}
\end{pmatrix}.
$$

Consider a unitary representation $\rho$ of $\G$
in a Hilbert space $H$. 
For any $n\ge m$ we have an operator
\begin{equation}
\wt\rho_{nn}\left(X\left(\begin{pmatrix}
b&0\\0&1_{n-m}
\end{pmatrix}, \begin{pmatrix}
c&0\\0&1_{n-m}
\end{pmatrix}\right)
\right):H_n\to H_n.
\label{eq:b1c1}
\end{equation}

We claim that these operators  as operators $H\to H$
depend only on $L$, $M$ and not on $n$. Precisely, we have the following statement.

\begin{proposition}
\label{pr:coherence}
	{\rm a)} Let $n\ge m$.
Then a block matrix structure of the operator {\rm(\ref{eq:b1c1})} with respect 
to the orthogonal decomposition
$H_n=H_m\oplus (H_n\ominus H_m)$
is
\begin{equation}
\wt\rho_{nn}\left( X\left(\begin{pmatrix}
b&0\\0&1_{n-m}
\end{pmatrix}, \begin{pmatrix}
c&0\\0&1_{n-m}
\end{pmatrix}\right)\right)=
\begin{pmatrix}
\wt\rho_{mm}(X(b,c))&0\\0&0
\end{pmatrix}.
\label{eq:XXX}
\end{equation}

{\rm b)}	For any $L$, $M\subset \frl^m$  we have a well-defined operator
	$\wt \rho(\cX[L,M])$
	in $H$, which  sends $H_m$ to $H_m$
	as $\wt \rho_{mm}(\cX[L,M])$ and is zero on the orthocomplement $H\ominus H_m$.
\end{proposition}

{\sc Proof.} According Corollary \ref{cor:iota}, the right hand side
of (\ref{eq:XXX})
is $\wt\rho_{nn}(X_m^\square(b,c))$. By Proposition \ref{pr:coincidence},
this operator coincides with $\wt\rho_{nn}(X_m^\bigcirc(b,c))$.
\hfill $\square$

\sm

{\bf\punct Proof of Proposition \ref{pr:coincidence}.%
\label{ss:skleika}}

\begin{lemma}
\label{l:fixed}
	Let $\frg\in \red(\Gamma(m))$ be an idempotent. Let
	$g$ be a representative of $\frg$ in $\G_\fin$. Then
	for any unitary representation $\rho$ of $\G$ in a Hilbert space $H$
	the image of the orthogonal projector $\wt\rho_{mm}(\frg)$ coincides with the space of fixed points
	of the subgroup in $\G$ generated by $\G(m)$ and $g$.
\end{lemma}

{\sc Proof.} Let $v\in\im\wt\rho_{mm}(\frg)$, i.e.,
$$
P_m\rho(g)P_mv=v.
$$
This happens if and only if $P_mv=v$, $\rho(g)v=v$. The condition $P_mv=v$
means that $\rho(h)v=v$ for all $h\in \G(m)$.
\hfill $\square$

\sm 

Therefore, it is sufficient to show that the group  generated by $\G(m)$
and $X^\bigcirc(b,c)$ coincides with the group generated by $\G(m)$
and $X^\square(b,c)$.

\begin{lemma}
	\label{l:generated}
The group  generated by the subgroup $\G(\beta)$
and the matrix
$$
X(1,1)=
\left(
\begin{array}{c|ccc}
1_\beta&1_\beta&0\\
\hline
0&1_\beta&0\\
1_\beta&0&1_{\beta\star}
\end{array}
\right)
$$	
coincides with $\G$.
\end{lemma}

{\sc Proof.}
Denote by $G$ the group generated by $X(1,1)$ and $\G(\beta)$.
Conjugating $X(1,1)$ by block diagonal matrices
we can get any matrix of the form $X(A,B)$ with nondegenerate $A$, $B$.
Multiplying such matrices we observe that elements of the form
$X(A_1+A_2,B_1+B_2)$ are contained in $G$.  In particular,  $X(0,2)\in G$. Since $p\ne 2$,
conjugating $X(0,2)$ by a block scalar matrix we come to $X(0,1)\in G$.
In the same way $X(1,0)\in G$. Now the statement became more-or-less obvious.
\hfill $\square$

\begin{lemma}
	\label{l:generated2}
	The group generated by $\G(\beta)$ and the matrix
	\begin{equation}
	\left(\begin{array}{c|c}
	0&1_\beta\\\hline 1_\beta& 0_\star
	\end{array}\right)
	\label{eq:tran}
	\end{equation}
	coincides with $\G$.  
\end{lemma}

{\sc Proof.} Denote this group by $G$. Denote $S_\infty(\beta):=S_\infty\cap \G(\beta)$.
 Multiplying
the matrix (\ref{eq:tran}) from the left and right by elements of
$S_\infty(\beta)$ we can get an arbitrary matrix of the form
$	\left(\begin{array}{c|c}
0&\sigma_1\\\hline \sigma_2& 0_\star
\end{array}\right)$ with  $\sigma_1$, $\sigma_2\in S_\beta$. Multiplying two matrices of this type
we can get any matrix $	\left(\begin{array}{c|c}
\sigma&1\\\hline 0& 1_\star
\end{array}\right)$, where $\sigma\in S_\beta$.
Therefore our group contains the subgroup $S_\beta\times S_\infty(\beta)$,
which is maximal in $S_\infty$. Therefore $G\supset S_\infty$.
 But  $S_\infty$ and $\G(\beta)$ generate $\G$, see \cite{Ner-p}, Lemma 3.6.
 \hfill $\square$

\sm

{\sc Proof of Proposition \ref{pr:coincidence}.} Denote by

\sm

--- $G^\bigcirc$  the group generated by $\G(m)$
and $X^\bigcirc_\alpha(b,c)$;

\sm 

--- $G^\square$   the group generated by $\G(m)$
and $X^\square_\alpha(b,c)$;

\sm  

--- $G$  the group generated by $\G(\alpha)$ and the matrix $X_\diamond(b,c)$
defined by
$$
X_\diamond(b,c):=
\left(\begin{array}{cc|cccc}
1_{\alpha}&0&b&0\\
0&1_{m-\alpha}&0&0\\
\hline
0&0&1&0\\
c&0&0&1_\star
\end{array}\right).
$$

\sm 

Obviously, $G\supset G^\bigcirc$, $G\supset G_\square$.
Let us verify the opposite inclusions.

\sm

{\it The inclusion $G^\bigcirc\supset G$.}
Clearly $X_\alpha(-b,- c)\in G^\bigcirc$. Therefore
$G^\bigcirc$ contains
$$
X_\alpha(b,c)X_\alpha(-b,-c)=
X\left(\begin{pmatrix}
0&0\\0&2
\end{pmatrix},
\begin{pmatrix}
0&0\\0&2
\end{pmatrix}
\right)
\sim
X\left(\begin{pmatrix}
0&0\\0&1
\end{pmatrix},
\begin{pmatrix}
0&0\\0&1
\end{pmatrix}
\right)=:Y.
$$
By Lemma \ref{l:generated}, the group generated by $Y$ and $\G(m)$
is $\G(\alpha)$. On the other hand, $Y^{-1}X^\bigcirc(b,c)\sim X_\diamond(b,c)$.

\sm 

{\it The inclusion $G^\square\supset G$.}
We have 
$$
X^\square_\alpha(b,c)^2\sim X_\diamond(2b,2c)\sim X_\diamond(b,c).
$$
Next, $X_\diamond(b,c)^{-1} X^\square_\alpha(b,c)\sim X_\alpha(0,0)$
and we refer to Lemma \ref{l:generated2}.

\sm

Thus, $G^\bigcirc=G^\square$. By Lemma \ref{l:fixed},
for any unitary representation $\rho$ of $\G$ we have
$$
\wt\rho_{mm}\bigl(X^\bigcirc(b,c))\bigr)=\wt\rho_{mm}\bigl(X^\square(b,c))\bigr)
$$
and this completes the proof of Proposition \ref{pr:coincidence}.
\hfill $\square$

\section{The semigroup $\red_m(\Gamma(m))$%
\label{s:red-m}}

\COUNTERS

{\bf \punct Structure of the semigroup $\red_m(\Gamma(m))$.}
Denote by $\0$ the minimal idempotent of the semigroup  $\red_m(\Gamma(m))$.

\begin{proposition}
	\label{pr:structure}
Any element $\ne \0$ in $\red_m(\Gamma(m))$ has a representative of a form
$a X(b,c)$, where $a\in \GL(m)$. 
\end{proposition}

The proof occupies the rest of the section.
As a byproduct of Lemma \ref{l:canonical-forms}
we will get the following statement. 

\begin{lemma}
	\label{l:idempotent-canonical}
Any idempotent $[X(b,c)]$ by a conjugation by $a\in \GL(m)$
 can be reduced to a form
$$
\left[X\left(\begin{pmatrix}0&\beta\\1&0 \end{pmatrix}, \begin{pmatrix}0&1\\\gamma&0 \end{pmatrix} \right) \right],
$$
where $\gamma\beta=0(\md p)$, $\beta\gamma=0(\md p)$.
\end{lemma}

\sm

{\bf\punct Proof of Proposition \ref{pr:structure}.}
{\sc Step 1.}

\begin{lemma}
	\label{l:canonical-forms}
	{\rm a)} Let $B$ be an $m\times N$ matrix over $\ZZ$,
	$C$ an $N\times m$ matrix. Then   transformations
	$$
B\mapsto u^{-1} B v,\qquad C\mapsto v^{-1} Cu
$$
allow to reduce them to the form
\begin{equation}
\wt B=\begin{pmatrix}
0&b_{12}\\
b_{21}&b_{22}
\end{pmatrix},
\qquad 
\wt C=
\begin{pmatrix}
0&c_{12}\\
c_{21}&c_{22}
\end{pmatrix}
\label{eq:BC}
,\end{equation}
where $b_{12}$, $c_{21}$ are square nondegenerate matrices of 
the same size, products $b_{21}c_{12}$,  $c_{12}b_{21}$
are nilpotent and $b_{22}=0(\md p)$, $c_{22}=0(\md p)$.

\sm 

{\rm b)} The transformations 
$$
B\mapsto u^{-1} B v,\qquad C\mapsto w^{-1} Cu,
$$
where $u$, $v$, $w$ are invertible, allow to reduce a pair
 $(B,C)$ to the form
 $$
 \wt B=\begin{pmatrix}
 0&b_{12}\\
 1&0
 \end{pmatrix},
 \qquad 
 \wt C=
 \begin{pmatrix}
 0&1\\
 c_{21}&0
 \end{pmatrix}
 ,$$
 where $c_{21}b_{12}=0(\md p)$, $b_{12}c_{21}=0(\md p)$.
\end{lemma}

{\sc Proof.} a) Reduce our matrices modulo $p$.
A canonical form of a pair of counter 
operators 
 $P:\F_p^m\to \F_p^N$ and
$Q:\F_p^N\to \F_p^m$ is a standard problem of linear algebra,
see, e.g., \cite{DP}, \cite{HS}. In particular, such operators in some bases admit
block decompositions $P=\begin{pmatrix}P_r&0\\0&P_n\end{pmatrix}$, 
$Q=\begin{pmatrix}Q_r&0\\0&Q_n\end{pmatrix}$, where $P_rQ_r$, $Q_rP_r$
are nondegenerate and $P_nQ_n$, $Q_nP_n$ are nilpotent.

Thus the matrices $B$, $C$ can be reduced to the form
$$
 B'=\begin{pmatrix}
b_{11}&b_{12}\\
b_{21}&b_{22}
\end{pmatrix},
\qquad 
 C'=
\begin{pmatrix}
c_{11}&c_{12}\\
c_{21}&c_{22}
\end{pmatrix},
$$
where 

\sm

1) $b_{21}$, $c_{12}$ are invertible matrices of the same size;

\sm

2) products $b_{12}c_{21}$, $c_{21}b_{12}$ are nilpotent;

\sm

3) the matrices $b_{11}$, $b_{22}$, $c_{11}$, $c_{22}$ reduced $(\md p)$ are zero.

\sm

Set
$$
u_1:=\begin{pmatrix}
1&b_{11} b_{12}^{-1}\\
0&1
\end{pmatrix},
$$
notice that  $u_1$ ($\md p)$ is 1.
We pass to new matrices
$$
B''= u_1^{-1} B',\qquad C''=C''u_1.
$$
For new $B$ the block $b_{11}=0$, other properties 1)--3) of matrices $B$, $C$
are preserved. Next, we take a unique  matrix of the form
$u_2=\begin{pmatrix}
1&0\\
*&1
\end{pmatrix}$
such that $C'' u_2$ has zero block $c_{12}$.
On the other hand the block $b_{11}$ of $u_2^{-1} B''$ is
zero.
 We come to a desired
form.

\sm

b) We apply  statement a) and reduce $(B,C)$ to the form
(\ref{eq:BC}). Next, we multiply $\wt B$ from right by 
$\begin{pmatrix}
b_{21}&\\&1
\end{pmatrix}^{-1}
$
and get 1 on the place of $b_{21}$. After this, we multiply new $B$ from right by
$\begin{pmatrix}
1&-b_{22}\\0&1
\end{pmatrix}$ and kill $b_{22}$. Finally, we repeat the same transformations with
$\wt C$.

Now the problem is reduced to the same question for a pair $b_{12}$, $c_{21}$.
If $c_{21}b_{12}\ne 0(\md p)$, then we choose an invertible matrix
$U$ such that $b_{12}Uc_{21}$ is not nilpotent and again repeat a). Etc.
\hfill $\square$

\sm

{\sc Step 2.}

\begin{lemma}
Let $[g]\in \Gamma(m)$ have the form
$$
[g]=\left[\begin{array}{c|c}
1&b\\
\hline
c&1_\star
\end{array}\right]_{mm}
$$
and $[[g]]\not\thickapprox_m\0$.
Then $bc$ and $cb$ are nilpotent.	
\end{lemma}

{\sc Proof.} We apply the previous lemma and represent
$[g]$ as
$$
[g]=
\left[\begin{array}{cc|cc}
1_\alpha&0&0&b_{12}\\
0&1_{m-\alpha}&b_{21}&b_{22}\\
\hline
0&c_{12}&1_{m-\alpha}&0\\
c_{21}&c_{22}&0&1_\star
\end{array}\right]_{mm}.
$$
Set 
$$[h_m^\alpha]:=
\left[\begin{array}{cc|cc}
1_\alpha&0&0&0\\
0&1_{m-\alpha}&1_{m-\alpha}&0\\
\hline
0&0&1_{m-\alpha}&0\\
0&1_{m-\alpha}&0&1_{m-\alpha\star}
\end{array}
\right].
$$
Let us show that 
\begin{equation}
[g]\circ [h_m^\alpha]\sim [g].
\label{eq:ghtheta}
\end{equation}
Indeed,
\begin{multline}
[g]\circ [h_m^\alpha]=
\left[
\begin{array}{cc|cccc}
1_\alpha&0&0&b_{12}&0&0\\
0&1_{m-\alpha}&b_{21}&b_{22}&1&0\\
\hline
0&c_{12}&1_{m-\alpha}&0&c_{12}&0\\
c_{21}&c_{22}&0&1&c_{22}&0\\
0&0&0&0&1&0\\
0&1&0&0&0&1_\star
\end{array}
\right]_{mm}
\sim\\\sim
\left[
\begin{array}{cc|cccc}
1_\alpha&0&0&b_{12}&0&0\\
0&1_{m-\alpha}&b_{21}&b_{22}&\boxed{1_{m-\alpha}}&0
\vphantom{\biggr|}
\\
\hline
0&c_{12}&1_{m-\alpha}&0&0&0\\
c_{21}&c_{22}&0&1&0&0\\
0&0&0&0&1&0\\
0&\boxed{1_{m-\alpha}}&0&0&0&1_\star
\end{array}
\right]_{mm}=: r,
\label{eq:h-circ-theta}
\end{multline}
to establish the equivalence we multiply
 $[g]\circ [h_m^\alpha]$ from the left
 by
 $$\left(
 \begin{array}{c|cccc}
 1_m&&&&\\
 \hline
 &1&0&-c_{12}&\\
 &0&1&-c_{22}&\\
 &0&0&1&\\
 &&&&1_\star
 \end{array}
 \right).
 $$
 Next, denote
 $$
 v_1:=
 \left(
 \begin{array}{c|cccc}
 1_m&&&&\\
 \hline
 &1&0&-b_{21}^{-1}&\\
 &&1&0&\\
 &&&1&\\
 &&&&1_\star
 \end{array}
 \right),
 \qquad 
 v_2:=
 \left(
 \begin{array}{c|cccc}
 1_m&&&&\\
 \hline
 &1&&&\\
 &&1&&\\
 &&0&1&\\
 &&-c_{12}^{-1}&0&1_\star
 \end{array}
 \right).
 $$
 We have
 $$
 [r]\sim [v_2v_1^{-1}r v_1v_2^{-1}],
 $$
the latter matrix is obtained from $r$, see (\ref{eq:h-circ-theta}),
by removing two boxed blocks $\boxed{1_{m-\alpha}}$, all other
blocks are the same. 
Thus $[r]\sim [g]$, i.e., we established (\ref{eq:ghtheta}).

\sm

Suppose that $\alpha\ne m$. Then by Proposition \ref{pr:coincidence},
$$
[g]\sim [g]\circ [h_m^{\alpha}]\thickapprox
[g]\circ \Theta_{[m]}^\alpha .
$$
But $\Theta_{[m]}^\alpha \thickapprox_m\0$, therefore $[[g]]\approxeq_m \0$.
\hfill $\square$

\sm

{\sc Step 3.} Thus it is sufficient to prove
Proposition \ref{pr:structure} for $[g]$ having the form
$$
[g]=\left[\begin{array}{c|c}
1&b\\
\hline
c&1_\star
\end{array}\right]_{mm},\qquad \text{where $bc$, $cb$ are nilpotent.}
$$

\begin{lemma}
	\label{l:gg-1}
Let $[g]=\left[\begin{array}{c|c}
1&b\\
\hline
c&1_\star 
\end{array}\right]_{mm}$ be invertible%
\footnote{This is equivalent to invertibility of $(1-bc)^{-1}$
	or invertibility of $(1-cb)^{-1}$. Here we do not need a nilpotency
of $bc$.}.	 Then 
	$$	[g^{-1}]\circ [g]\equiv X(b,c).$$
\end{lemma}

{\sc Proof.}
By (\ref{eq:xxx}), 
$$
[g]\circ ([g^{-1}]\circ [g])\equiv[g], 
\qquad\text{and $[g^{-1}]\circ [g]$ is an idempotent.}
$$
We have (see, e.g., \cite{Gan}, Sect. 2.5)
$$
[g^{-1}]=
\left[
\begin{array}{c|c}
(1-bc)^{-1}&-(1-bc)^{-1}b\\
\hline
-c(1-bc)^{-1}&(1-cb)^{-1}_\star
\end{array}
\right]_{mm}
\sim
\left[
\begin{array}{c|c}
(1-bc)^{-1}&b\\
\hline
c(1-bc)^{-1}&1_\star
\end{array}
\right]_{mm}.
$$
We also keep in mind the identity
\begin{equation}
c(1-bc)^{-1}=(1-cb)^{-1}c,
\label{eq:1-cb}
\end{equation}  
to establish it, we multiply both sides from the left by
$(1-cb)$ and from the right by $(1-bc)$.

Next,
$$
[g^{-1}]\circ [g]=
\left[
\begin{array}{c|cc}
(1-bc)^{-1}&b&(1-bc)^{-1}b\\
\hline
c(1-bc)^{-1}&1&c(1-bc)^{-1}b\\
c&0&1_\star
\end{array}\right]_{mm}
\!\!\!\!
\sim
\left[
\begin{array}{c|cc}
(1-bc)^{-1}&b&b\\
\hline
c&1&0\\
c(1-bc)^{-1}&0&1_\star
\end{array}\right]_{mm}
.
$$
This matrix defines an idempotent in
$\inv(\Gamma^\circ(m))$. We must verify the following statement:

\begin{lemma}
	\label{l:identification}
	Under our conditions,
	$$[g^{-1}]\circ [g]\equiv X(b,c).$$
\end{lemma}

{\sc Proof.}
 By Corollary \ref{cor:}
we can identify an idempotent in $\inv(\Gamma^\circ)$ evaluating its image in $\Mor_\cL(m)$.
So we get
$$
[g]\circ [g^{-1}]\equiv[[X(B,C)]],
$$
where
$$
B:=\begin{pmatrix}
b&b
\end{pmatrix},\qquad 
C:=\begin{pmatrix}
c\\(1-cb)^{-1}c
\end{pmatrix}.
$$
We have $\ker B=\ker b$, $\ker C^t=\ker c^t$,
therefore by Lemma \ref{l:idempotents-2} we have
$[[X(B,C)]]\equiv[[X(b,c)]]$.
\hfill$\square$

\begin{corollary}
Let 
$$[g]=\left[\begin{array}{c|c}
	1&b\\
	\hline
	c&1_\star
\end{array}\right]_{mm},\qquad
[g']=\left[\begin{array}{c|c}
1&bu\\
\hline
c&1_\star
\end{array}\right]_{mm}	
$$
be invertible and $u$ also be invertible.
Then 
$$
[g^{-1}]\circ [g]\equiv[(g')^{-1}]\circ [g'].
$$
\end{corollary}

{\sc Proof.} Indeed, $\ker bu=\ker b$. So both sides are
$[[X(b,c)]]$.
\hfill $\square$

\sm

{\sc Step 4.}

\begin{lemma}
	\label{l:2}
Let $[g]=\left[\begin{array}{c|c}
1&b\\
\hline
c&1_\star
\end{array}\right]_{mm}$, let $bc$ and $cb$ be nilpotent.
Then there exists $u$
having the form
\begin{equation}
u=-\frac12+\sum_{j>0} \frac{\sigma_j}{2^{n_j}} (cb)^j,
\qquad \text{where $\sigma_j\in \Z$, $n_j\in\Z_+$,}
\label{eq:sum}
\end{equation}
such that 
\begin{multline}
\left(
\left[\left(\begin{array}{c|c}
1&bu\\
\hline
c&1_\star
\end{array}\right)^{-1}\right]_{mm}\circ
\left[\begin{array}{c|c}
1&bu\\
\hline
c&1_\star
\end{array}\right]_{mm}
\right)\circ [g^{-1}]\equiv 
\\\equiv
\left[
\begin{array}{c|cc}
(1-buc)^{-1}(1-bc)^{-1}&b&0\\
\hline
0&1&0\\
c&0&1_\star
\end{array}\right]_{mm}
.
\end{multline}
\end{lemma}

{\sc Proof.} The product is 
\begin{multline*}
\left[
\begin{array}{c|ccc}
(1-buc)^{-1}(1-bc)^{-1}&bu&bu&(1-buc)^{-1}b\\
\hline
c(1-bc)^{-1}&1&0&cb\\
c(1-buc)^{-1}&0&1&c(1-buc)^{-1}\\
c(1-bc)^{-1}&0&0&1_\star
\end{array}\right]_{mm}
\sim\\\sim
\left[
\begin{array}{c|ccc}
(1-buc)^{-1}(1-bc)^{-1}&bu&bu&(1-buc)^{-1}b\\
\hline
c&1&0&0\\
c(1-buc)^{-1}&0&1&0\\
c(1-bc)^{-1}&0&0&1_\star
\end{array}\right]_{mm}=:
\left[
\begin{array}{c|c}
A&br\\
\hline
qc&1_\star
\end{array}\right]_{mm},
\end{multline*}
here 
$$
r:=\begin{pmatrix}
u&u&(1-cbu)^{-1}
\end{pmatrix},\qquad
q:=
\begin{pmatrix}
1\\
(1-cbu)^{-1}\\
(1-cb)^{-1}
\end{pmatrix}.
$$
We claim that {\it there exists a unique $u$ such that
$rq=0$}. A straightforward calculation shows that
$$
rq=2u-ucbu+(1-cb)^{-1}.
$$
Since $cb$ is nilpotent, we can write the equation $rq=0$
as
$$
2u+1=ucbu-\sum_{j>0}(cb)^j,
$$
the sum actually is finite. Clearly
we can find a solution in the form
$u=-1/2+\sum_{j>0} s_j(cb)^j$, where $s_j$ are dyadic rationals,
for coefficients $s_j$ we have a system of recurrent equations.
This $u$ is invertible (since we can write a finite series for
$u^{-1}$).

Next, we must show that the matrix
$\left(\begin{array}{c|c}
1&bu\\
\hline
c&1_\star
\end{array}\right)$ is invertible. Indeed, this is equivalent
to existence of
$
(1-cbu)^{-1}
$
and this is clear since by (\ref{eq:sum}) $cbu$ is nilpotent.

\sm

Next we wish to simplify
the matrix
$\left(\begin{array}{c|c}
	A&br\\
	\hline
	qc&1_\star
\end{array}\right)$ by conjugations by matrices of the form
$
\begin{pmatrix}
1&0\\0&D_\star
\end{pmatrix}
$. In fact, we have transformations
$$
r\mapsto r'=pD^{-1},\qquad q\mapsto q'=Dq.
$$
For such transformations we have
$r'q'=rq$. Set
$$
D=\begin{pmatrix}
1&1&u^{-1}(1-cbu)^{-1}\\
0&1&0\\
0&0&1_\star
\end{pmatrix}.
$$
Then $r'=\begin{pmatrix}
u&0&0
\end{pmatrix}$.  But $u$ is invertible and $r'q'=0$. Therefore
$q'$ has the form
$\begin{pmatrix}
0\\ *\\ *
\end{pmatrix}$, on the other hand multiplication $q\mapsto Dq$
does not change the second and third elements of the column $q$.
Thus we came to the matrix
$$R:=
\left[
\begin{array}{c|ccc}
(1-buc)^{-1}(1-bc)^{-1}&b\,\boxed{u}&0&0
\vphantom{\biggr|}
\\
\hline
0&1&0&0\\
\boxed{(1-cbu)^{-1}}\,c&0&1&0
\vphantom{\Bigr|}
\\
\boxed{(1-cb)^{-1}}\,c&0&0&1_\star
\end{array}\right]_{mm}.
$$
Consider the following matrices:
$$
S:=
\left[\begin{array}{c|ccc}
1&&&
\\
\hline
&u&&\\
&&1-cbu&\\
&&&(1-cb)_\star
\end{array}\right],\qquad\qquad
\\
T:=
\left[\begin{array}{c|ccc}
1&&&
\\
\hline
&1&&\\
&&1&\\
&&-1&1_\star
\end{array}\right].
$$
The conjugation $R\mapsto TRT^{-1}$ kills boxed elements of $R$.
The conjugation $R\mapsto STRT^{-1}S^{-1}$
reduces the matrix to the desired form.
\hfill $\square$

\sm

{\sc Proof of Proposition \ref{pr:structure}.}
Thus we have
$$
\left[\left(\begin{array}{c|c}
1&b\\
\hline
c&1_\star
\end{array}\right)^{-1}\right]_{mm}
\equiv
(1-buc)^{-1}(1-bc)^{-1}
\cdot [[X\bigl((1-bc)(1-buc)b,c\bigr)]].
$$
The second factor is
$$
[[X\bigl(b(1-cb)(1-cbu),c\bigr)]]\equiv [[X(b,c)]].
$$
Passing to adjoint elements we get
\begin{multline*}
\left[\begin{array}{c|c}
1&b\\
\hline
c&1
\end{array}\right]_{mm}\equiv[[X(b,c)]]\cdot
(1-bc)(1-buc)
\equiv\\\equiv
 (1-bc)(1-buc)
\cdot[[X((1-buc)^{-1}(1-bc)^{-1}b,c (1-bc)(1-buc) )]]
 \equiv\\\equiv
 (1-bc)(1-buc)
\cdot [[X(b,c)]].
\end{multline*}
It remains to notice that
$$
\left[\begin{array}{c|c}
a&b\\
\hline c&1
\end{array}\right]
=a\cdot
\left[\begin{array}{c|c}
1&a^{-1}b\\
\hline c&1
\end{array}\right].
$$

{\bf \punct Proof of Lemma \ref{l:idempotent-canonical}.}
We refer to Lemma \ref{l:canonical-forms}.

\section{The groups $\G^\bullet[L;M]$%
\label{s:bullet}}

\COUNTERS

In this section we examine subgroups $\G^\circ[L;M]$,
$\G^\bullet[L;M]\subset \G$ defined in Subsection \ref{ss:statement}.
We prove that $\G^\bullet[L;M]$ is well-defined (Lemma \ref{l:G-bullet-well},
shows that it is generated by $\G(m)$ and the idempotent $X(b,c)$ (Proposition
\ref{pr:bullet}). 
Also we prove that it
is a minimal subgroup of finite index
in $\G[L;M]$ (equivalently, $\G^\bullet[L;M]$ has no subgroups of finite index,
Proposition \ref{pr:no-finite}).

\sm

{\bf\punct Several remarks on submodules in $\frl^k$.}

\begin{lemma}
	\label{l:can-basis}
	 Let $L\subset \frl^k$ be a submodule. Then there exists a basis
	 $e_j\in  \frl^k$ such that $M:=\oplus p^{s_j} \ZZ e_j$.
	 The collection $s_1$, $s_2$, \dots is a unique $\GL(m)$-invariant of a submodule $L$. 
\end{lemma}

This is equivalent to a classification of sublattices in $(\O_p)^k$ under the action
of $\GL(k,\O_p)$ or equivalently to a classification of  pairs of  lattices
in $\Q_p^k$ under $\GL(k,\Q_p)$, the latter question
is standard, see, e.g., \cite{Wei}, Theorem I.2.2.
\hfill $\square$

\begin{corollary}
	Any submodule $L\subset \frl^k$ is a kernel of
	some  endomorphism $\frl^k\to \frl^k$.
\end{corollary}

Indeed, we pass to a canonical basis $e_j$ as in the lemma and consider the map  sending
$e_j$ to $p^{\mu-s_j}e_j$.
\hfill $\square$



\begin{lemma}
	\label{l:52}
	{\rm a)}
	Let $L$ be a submodule in $\frl^m$. Let 
	$b$, $b':\frl^m\to \frl^N$ be  morphisms of modules
	such that $L=\ker b=\ker b'$. Then there is a transformation
	$u\in \GL(N)$ such that $b'=bu$.
	
	\sm
	
	{\rm b)} Let $\ker b=L$, $\ker b'=L'\supset L$. Then there is an endomorphism
	$u:\frl^N\to \frl^N$ such that $b'=bu$.
\end{lemma}

{\sc Proof.} a) The modules $\im b\simeq \im b'\simeq \frl^m/L$
are isomorphic. By the previous lemma there is an automorphism
of $\frl^N$ identifying these submodules.

\sm

b)  $L$ is a submodule of $L'$, therefore $\im b'$ is a quotient module
of $\im b$. Therefore there is a projection map $\pi:\im b\to\im b'$, orders of
elements do not increase under this map. 
By Lemma \ref{l:can-basis} we have a basis $e_j\in \frl^N$
such that $p^{s_j} e_j$, where $j=1$, \dots, $m$, is the system of generators of  
$\im b$. Choose arbitrary vectors $v_j$ such that 
$p^{s_j}v_j=\pi(p^{s_j} e_j)$ and consider the map sending $e_j$ to $v_j$.
\hfill $\square$.

\sm

{\bf\punct The group $\G^\bullet$.}
Here we show that $\G^\bullet[L;M]$ is a group,
and its definition does not depend on the choice of matrices $b$, $c$.

\begin{lemma}
	{\rm a)} Fix a  matrix $B$ of size $l\times N$. Then the set of invertible matrices $g$
	of the form $1-BS$, where $S$ ranges in the set of $N\times l$ matrices,
	is a group.
	
	\sm 
	
	{\rm b)} Fix matrices $B$, $C$ of sizes $l\times N$ and $N\times l$ 
	respectively. Then the set of invertible matrices $g$ of the form
	$g=1-BuC$ is a group.
\end{lemma}

{\sc Proof.} Clearly, both sets are closed with respect to
multiplication. We must show that $g^{-1}$ satisfies the same property.
In the first case,
$$
1-g^{-1}=1-(1-BS)^{-1}=-BS(1-BS)^{-1}.
$$
 In the second case,
$$
1-g^{-1}=1-(1-BuC)^{-1}=-BuC(1-BuC)^{-1}=
-Bu(1-CBu)^{-1}C.
\qquad \square
$$

\begin{lemma}
	\label{l:G-bullet-well}
	 Fix matrices $b$, $c$ of sizes $m\times N$ and $N\times m$
	respectively.
	
	\sm 
	
	{\rm a)}
	The set of invertible matrices $g=\begin{pmatrix} a&bv\\wc&z
	\end{pmatrix}$ such that the block '$a$' admits representations
	$a=1-bS$, $a=1-Tc$ is a group
	
	\sm
	
	{\rm b)} 
	 The set $\G^\bullet[L;M]$, i.e., the set of all invertible matrices
	 of the form 
	$g=\begin{pmatrix} 1-buc&bv\\wc&z
	\end{pmatrix}$,
	is a group.
\end{lemma}

{\sc Proof.}  In the first case
we write
\begin{multline*}
g=\begin{pmatrix}1-bS&bv\\wc&z\end{pmatrix}=
\begin{pmatrix}
1&\\&1
\end{pmatrix}-
\begin{pmatrix}
bS&-bv\\
-wc&1-z
\end{pmatrix}=\\=
\begin{pmatrix}
1&\\&1
\end{pmatrix}-
\begin{pmatrix}
b&\\&1
\end{pmatrix}
\begin{pmatrix}
S&-v\\
-wc&1-z
\end{pmatrix},
\end{multline*}
and  reduce the statement to the previous lemma.

In the second case we write
$$
g=
\begin{pmatrix}
1&\\&1
\end{pmatrix}-
\begin{pmatrix}
b&\\&1
\end{pmatrix}
\begin{pmatrix}
u&-v\\
-w&1-z
\end{pmatrix}
\begin{pmatrix}
c&\\&1
\end{pmatrix},
$$
and again we apply the previous lemma.
\hfill $\square$

\sm


{\bf \punct The group $\G^\circ[L;M]$.}

\sm 

{\sc Proof of Lemma \ref{l:G-circ}.}
Let $g=\begin{pmatrix}
\alpha&\beta\\\gamma&\delta_\star
\end{pmatrix}\in \G^\circ[L,M]$, i.e., $g$ fix pointwise $L\subset \frl^m$
and $g^t$ fix pointwise of $M\subset \frl^m$.
 Then 
$L\subset \ker \beta$ and by Lemma \ref{l:52}.b we have $\beta=b v$ for some matrix $v$. 
Also $L\subset \ker (1-\alpha)$ and therefore $\alpha=1-b S$ for some $S$.
\hfill $\square$

\sm 

{\bf\punct Changes of coordinates.}
\begin{lemma}
	Let $L$, $M\subset \frl^m$. 
	Let $a\in\GL(m)$. Then
$$ 
a \G^\circ[L;M] a^{-1}=
\G^\circ[aL,a^{t-1}M], \qquad
a \G^\bullet[L;M] a^{-1}=
\G^\bullet[aL,a^{t-1}M]
$$	
\end{lemma}

The first statement is an immediate consequence of the definition,
the second is straightforward. 
\hfill $\square$

\sm

{\bf\punct Generators of $\G^\bullet[L;M]$.}
 Let $m$, $b$, $c$ be as in Subsection \ref{ss:statement},
 i.e., $L=\ker b$, $M=\ker c^t\subset \frl^m$. 

\begin{proposition}
\label{pr:G-bullet-generators}
	The group  $\G^\bullet[L;M]$ is generated by
	$\G(m)$ and the matrix $X(b,c)$.
\end{proposition}

{\sc Proof.} Consider the group $G$ generated by $\G(m)$ and  $X(b,c)$.
Clearly, $\G^\bullet[L,M]\supset G$. Let us prove the converse.

\sm

 1) Conjugating $X(b,c)$ by block diagonal matrices
$\in \G(m)$ we get 
arbitrary matrices of the form $X(bv,wc)$, where $v$, $w$ are invertible matrices.
Consider products
\begin{equation}
X(bv,wc)\,X(b'v,wc')=X((b+b')v,w(c+c')).
\label{eq:XX}
\end{equation}
We set $b=-b'$, for any matrix $\sigma$ we can find invertible matrices
$c$, $c'$ such that%
\footnote{It is sufficient to verify this statement
	for matrices over $\F_p$. Without loss of generality
we can assume that $\sigma$ is diagonal. For $p\ne 2$
any element of $\F_p$ is a sum of two nonzero elements,
where $\sigma$ can be represented as a sum of two diagonal matrices.}
$c+c'=\sigma$.
Thus $G$ contains all matrices of the form
\begin{equation}
\begin{pmatrix}
1_m&bv\\
&1_{m\star}
\end{pmatrix},\qquad
\begin{pmatrix}
1_m&\\
wc&1_{m\star}
\end{pmatrix},
\label{eq:tt}
\end{equation}
where $v$, $w$ are arbitrary matrices.

\sm 

2) In virtue of Lemma \ref{l:idempotent-canonical}, conjugating the matrices (\ref{eq:tt}) by elements of $\GL(m)$ and multiplying from the left and the right
by elements of $\G(m)$
 we can reduce the matrices (\ref{eq:tt}) to the forms
\begin{equation}
Y[\beta]:=
\begin{pmatrix}
1_{m-\alpha}&0&0&\beta \\
0&1_\alpha&\boxed{1_\alpha}&0\\
0&0&1_{\alpha}&0\\
0&0&0&1_{m-\alpha\star}
\end{pmatrix},
Z[\gamma]:=
\begin{pmatrix}
1_{m-\alpha}&0&0&0\\
0&1_\alpha&0&0\\
0&\boxed{1_\alpha}&1_{\alpha}&0\\
\gamma &0&0&1_{m-\alpha\star}
\end{pmatrix},
\label{eq:111}
\end{equation}
where $\gamma\beta=0(\md p)$, $\beta\gamma=0(\md p)$.
Multiplying  $Y[\beta]$ from right by elements of $\G(m+\alpha)$
we can get that any matrix $Y[\beta r]$ with invertible $r$.
The condition of invertibility of $r$ can be removed, because
$$
Y[\beta r_1]\, Y[\beta r_2]^{-1} \, Y[\beta r_3]
=Y[\beta (r_1-r_2+r_3)],
$$ 
and we can represent any matrix $r$ as a sum of 3 invertible matrices.
 
 In the same way we get that $G$ contains all elements 
 of the form $Z[q\gamma]$.

Take $r=0$, $q=0$. Then the matrices $Y[0]=Y[\beta\cdot 0]$,
 $Z[0]=Z[0\cdot \gamma]$ together with $\G(m)$ generate
the group $\G(m-\alpha)$.

Next, $G$ contains  matrices $Y[\beta]Y[0]^{-1}$
and $Z[\gamma]Z[0]^{-1}$. They are matrices
 of the form (\ref{eq:111}),
where boxed blocks are replaced by zeroes. 

 Therefore our problem is reduced to
a description of the subgroup generated by
$\G(m-\alpha)$ and $X(\beta,\gamma)$.

{\it Thus, without loss of generality, we can assume that $\alpha=0$ and
$cb=0(\md p)$, $bc=0(\md p)$}.

\sm

3) Multiplying the matrices (\ref{eq:tt}), we get
$$
\begin{pmatrix}
1-bvwc&bv\\
wc&1_\star
\end{pmatrix}\in G \quad \text{for any $v$, $w$.}
$$
Since $cb=0(\md p)$, the $bvwc$ is nilpotent, and therefore
 $1-bvwc$ is invertible. We represent our matrix as
 \begin{multline*}
 \begin{pmatrix}
 1&0\\
 wc(1-bvwc)^{-1}&1_\star
 \end{pmatrix}
 \begin{pmatrix}
 1-bvwc&0\\0&1_\star
 \end{pmatrix}
 \times\\\times
 \begin{pmatrix}
 1&0\\0&1-wc(1-bvwc)^{-1}bv_\star
 \end{pmatrix}
 \begin{pmatrix}
 1&(1-bvwc)^{-1}bv\\
 0&1_\star
 \end{pmatrix}.
 \end{multline*}
 Since the whole product and three factors are contained in $G$,
 the fourth factor also is contained in $G$,
 $$
  \begin{pmatrix}
 1-bvwc&0\\0&1_\star
 \end{pmatrix}\in G
 $$
 for any $v$, $w$.
 
 \sm 
 
 4) Now consider an arbitrary element of $\G^\bullet[L;M]$,
 \begin{multline*}
 \begin{pmatrix}
 1-buc&bv\\
 wc&z_\star
 \end{pmatrix}=
 \begin{pmatrix}
 1&0\\wc(1-buc)^{-1}&1_\star
 \end{pmatrix}
 \begin{pmatrix}
 1-buc&0\\
 0&1_\star
 \end{pmatrix}
 \times\\\times
 \begin{pmatrix} 
 1&0\\0&z-wc(1-buc)^{-1}bv_\star
 \end{pmatrix}
 \begin{pmatrix}
 1&(1-buc)^{-1}bv\\0&1_\star
 \end{pmatrix}
 \end{multline*}
 All factors of the right hand side are contained in $G$,
 and therefore $\G^\bullet[L;M]$ is contained in $G$.
 \hfill $\square$

\begin{corollary}
The group $\G^\bullet[L;M]$ does not depend on a choice of $m$.	
\end{corollary}

{\sc Proof.} Let $L$, $M\subset \frl^m$, let
$L=\ker b$, $M=\ker c^t$. Let us regard $L$, $M$ as
submodules $L'$, $M'$ of $\frl^m\oplus \frl^k$. Then 
$$
\text{$L'=\ker b'$, $M'=\ker (c')t$, where  
	$b'=\begin{pmatrix}b&0\\0&1
\end{pmatrix}$, $c'=\begin{pmatrix}c&0\\0&1
\end{pmatrix}$}.
$$
Clearly the subgroup $G_m$ generated by $\G(m)$ and $X(b,c)$
and the subgroup $G_{m+k}$ generated by $\G(m+k)$ and $X(b',c')$
coincide. Formally, we must repeat the first two steps of the previous
proof.
\hfill $\square$

\sm
 
 {\bf \punct The quotient $\G^\circ/\G^\bullet$.}

 \begin{lemma}
 	A group $\G^\bullet[L;M]$ has  finite index in $\G^\circ[L;M]$.
 \end{lemma}

{\sc Proof.}
Without loss of generality we can assume that $cb=0(\md p)$,
$bc=0(\md p)$.
Denote by $A^\circ\subset \GL(m)$ the subgroup
consisting of matrices $a$ admitting  representations
$a=1-bS$, $a=1-Tc$. Notice that $1-a$ is a nilpotent,
since  $TcbS=0(\md p)$. Therefore $a$ is invertible.
Denote by $A^\bullet$ the subgroup consisting of elements
of the form $1-buc$. 

The subgroup $A^\bullet$ is normal in $A^\circ$.
Indeed, let $a\in A^\circ$, $a=1-bS$, $a^{-1}=1-Tc$. Then
$$
a(1-buc)a^{-1}=1-abuca^{-1}=
1-(1-bS)buc(1-Tc)=1-b(1-Sb)u(1-cT)c.
$$

Let $g=\begin{pmatrix}
a&bv\\wc & z
\end{pmatrix}\in \G^\circ[L;M]
$. Let us show that the map $g\mapsto a$ induces 
a homomorphism from $\G^\circ[L;M]\to A^\circ/A^\bullet$.
Indeed,
$$
g_1 g_2=\begin{pmatrix}
a_1&bv_1\\w_1 & z_1
\end{pmatrix}
\begin{pmatrix}
a_2&bv_2\\w_2c & z_2
\end{pmatrix}
=
\begin{pmatrix}
a_1 a_2+bv_1 w_2 c&*\\ * & *
\end{pmatrix}.
$$
In the left upper block we have
$$
a_1a_2(1+a_2^{-1}a_1^{-1}bv_1 w_2 c ).
$$
 We represent $a_1^{-1}=1-bS_1$, $a_2^{-1}=1-bS_2$ and get
$$
a_1a_2\bigl(1+(1-bS_2)(1-bS_1)bv_1 w_2 c \bigr)=
a_1a_2\bigl\{1+b(1-S_2b)(1-S_1b)v_1 w_2 c \bigr\}.
$$ 
The expression in the curly brackets is contained in $A^\bullet$.

Clearly, the kernel of the homomorphism is $\G^\bullet[L,M]$. Thus we have an isomorphism of quotient groups,
 $$
 \G^\circ[L;M]/\G^\bullet[L;M]\simeq A^\circ[L;M]/A^\bullet[L;M].
 $$
 The group in the right-hand side is finite.
 \hfill $\square$

\sm 

{\bf \punct Absence of subgroups of finite index.}

\begin{lemma}
	The group $\G$ has not proper open subgroups of finite index.
\end{lemma}

{\sc Proof.}
Let $P$ be a proper open subgroup. Then it contains
some group $\G(\nu)$. On the other hand $\G$ contains 
a complete infinite symmetric group $S_\infty$,
and $S_\infty$ has no subgroups of finite index.
Therefore  $P$ contains $S_\infty$. But the subgroup
in $\G$
generated by $\G(\nu)$ and $S_\infty$ is
the whole group $\G$, see \cite{Ner-p}, Lemma 3.6.
\hfill $\square$

\begin{proposition}
	\label{pr:no-finite}
The subgroup
$\G^\bullet[L;M]$ has no proper open subgroups of finite index.
\end{proposition}	  

{\sc Proof.} Let $Q$ be such subgroup. By the previous
lemma, $\G(m)$ has not open subgroups of finite
index, we have $Q\supset \G(m)$. Hence $Q$ contains a minimal
normal subgroup $R$ containing $\G(m)$. The quotient
$Q/R$ is generated by the image $\xi$ of $X(b,c)$,
therefore $Q/R$ is a cyclic group. But
$$X(b,c)^2=X(2b,2c)=
\left(
\begin{array}{c|cc}
1&&\\
\hline
&1/2&\\
&&2_\star
\end{array}
\right)
\left(
\begin{array}{c|cc}
1&b&0\\
\hline
0&1&0\\
c&0&1_\star
\end{array}
\right)
\left(
\begin{array}{c|cc}
1&&\\
\hline
&2&\\
&&1/2_\star
\end{array}
\right)
.$$
Since $p\ne 2$ the elements $X(b,c)^2$ and $X(b,c)$ have the same images
in $Q/R$. Therefore the image of $X[b,c]$ is 1.
\hfill $\square$

\begin{corollary}
	Any subgroup of finite index in $\G[L,M]$ contains $\G^\bullet[L,M]$.
\end{corollary}

\section{End of the proof%
\label{s:fin}}

\COUNTERS 

This section contains the end of the proof
of Theorem \ref{th:main}.
We know that all idempotents in semigroups $\red(\Gamma(n))$
have the form $\cX[L,M]$, see Corollary \ref{cor:all-idempotents},
 for different $n$ they can be identified in a natural way,
 see Proposition \ref{pr:coherence}. We also know that any non-zero
 element of $\red_m(\Gamma(m))$ is a product 
 of an invertible element and an idempotent $\cX[L,M]$, see Proposition \ref{pr:structure}.
 This  implies that all irreducible representations 
 of $\G$ are induced from representations $\tau$ of groups $\G[L;M]$.
 Proposition \ref{pr:no-finite} implies that  such $\tau$ must be trivial
 on $\G^\bullet[L;M]$.

\sm 

{\bf \punct A preliminary remark.}

\begin{lemma}
	\label{l:GNS}
Consider an irreducible $*$-representation $\sigma$
of the category $\cK$ in a sequence of Hilbert spaces
$H_j$. Let $\xi\in H_m$ be a nonzero vector.
Then the matrix element
$$c(\frg)=\la\sigma(\frg)\xi,\xi\ra_{H_m},\qquad
\text{where $\frg$ ranges in $\End_\cK(m)$,}$$
determines $\sigma$ up to equivalence.
\end{lemma}

This is a general statement on $*$-representations of categories
(and a  copy of a similar statement for unitary representations
of groups),
we give a proof for completeness.

\sm 

{\sc Proof.} 
For each  $\frg\in \Mor_\cK(m,\alpha)$
we define a vector 
$$\omega_{\frg}^\alpha=\sigma(\frg)\xi\in H_\alpha$$
Since $\sigma$ is irreducible, vectors $\omega_{ \frg}^\alpha$,
where  $\frg$ ranges in $\Mor_\cK(m,\alpha)$,
generate the  space $H_\alpha$.
Their inner products are determined by the function $c$:
$$
\la \omega_{ \frg_1}^\alpha,\omega_{ \frg_2}^\alpha\ra_{H_\alpha}=
\la \sigma(\frg_1)\xi, \sigma(\frg_2)\xi\ra_{H_\alpha}=
\la \sigma(\frg_2^*\circ\frg_1) \xi, \xi\ra_{H_m}=c(\frg_2^*\circ\frg_1).
$$
Next, let   $\frh\in \Mor_\cK(\alpha,\beta)$. Let $\frg$, $\frf$ range
respectively 
in $\Mor_\cK(m,\alpha)$, $\Mor_\cK(m,\beta)$. 
Then 
$$
\la \sigma(\frh)\omega_{ \frg},\omega_{\frf}\ra_{H_\beta}
=\la \sigma(\frh)\sigma(\frg)\xi, \sigma(\frf)\xi\ra_{H_\beta}
=\la \sigma(\frf^*\circ \frh\circ\frg)\xi, \xi\ra_{H_m}
=c(\frf^*\circ \frh\circ\frg).
$$ 
Clearly an operator $\sigma(\frh)$ is uniquely determined by such inner products.
\hfill $\square$

\sm 

{\bf \punct Representations of the semigroup $\red_m(\Gamma(m))$.}
Consider an irreducible representation of $\cK$ of height
$m$ and the corresponding representation $\lambda$ of the semigroup
$\End_\cK(m)$ in $H_m$. Recall that $\tau$ passes
through semigroup $\red_m(\Gamma(m))$. 
By Proposition \ref{pr:structure}, any
nonzero element
of the latter semigroup can be represented
as $a\cdot \cX[L,M]$, where $a\in \GL(m)$.
Denote
$$
\wh \G_n[L,M]=\GL(n)\cap \wh \G[L,M],\qquad 
\wh \G_\fin[L,M]=\G_\fin\cap \wh \G_m[L,M].
$$

The following lemma is a special case of general description 
of representations of finite inverse semigroups, see, e.g.,
\cite{GMS}. However, due to Proposition \ref{pr:structure}
our case is simpler than general inverse semigroups.
We show that the representation of $\GL(m)$ in
$H_m$ is induced from an irreducible representation of 
some subgroup $\wh \G_m[L,M]$ and idempotents
$\cX[N,K]$ act in the induced representation as multiplications
by indicator functions of certain sets. Precisely,

\begin{lemma}
	\label{l:description}
	Let $\cX[L,M]$ be the minimal idempotent 
	in $\red_m(\Gamma(m))$ such that $\lambda(\cX[L,M])\ne 0$.
	 Then there is an irreducible
	representation $\tau_m$ of $\wh \G_m[L,M]$ in a space $V$ such that
	$H_m$ can be identified with the space $\ell_2$ of $V$-valued 
	functions on the homogeneous space
	 $\wh \G_m[L,M]\setminus\GL(m)$ and
	
	\sm 
	
	{\rm 1)}  The group $\GL(m)$ acts by transformations
	of the form
	$$
	\lambda(g) f(x)= R(g,x) f(xg), 
	$$
	and for $q\in \wh \G_m[L,M]$ we have $R(p,x_0)=\tau_m(q)$ {\rm(}where $x_0$ 
	denotes the initial
	point of  $\wh \G_m[L,M]\setminus\GL(m)${\rm)}.
	
	\sm 
	
	{\rm 2)} The semigroup of idempotents acts by multiplications
	by indicator functions. Namely $\cX[K,N]$
	acts by multiplication by the function
	$$
	I_{K,N}(x_0 a)=\begin{cases}1,\quad\text{if
	  $K\supset aL$, $N\supset a^{t-1}M$};
  \\
  0,\quad\text{otherwise}.
	\end{cases}
	$$
\end{lemma}


\sm

{\sc Proof.} 
 Consider the image $V$ of the projector
$\lambda\bigl(\cX[L,M]\bigr)$. The idempotent
$\cX[L,M]$ commutes with $\wh \G_m[L,M]$. Indeed,
for $q\in \wh \G_m[L,M]$ we have
$$q\cdot\cX[L,M]\cdot q^{-1}=\cX[Lq, Mq^{t-1}]=\cX[L,M].$$
Therefore
the subspace $V$ is $\wh \G_m[L,M]$-invariant.
Denote by $\tau_m$ the representation of the group $\wh \G_m[L,M]$ in $V$.
We need the following lemma:

\begin{lemma}
	\label{l:induced}
For any $\frg\in \red_m(\Gamma(m))$ we have $\lambda(\frg) V=V$
or $\lambda(\frg) V\bot V$.
\end{lemma}

{\sc Proof of Lemma \ref{l:induced}.}
Let us apply an arbitrary element of $\red_m(\Gamma(m))$
to $v\in V$,
$$
\lambda\left(a\cdot \cX[K,N]\right)v=\lambda(a)\cdot \lambda\left(\cX[K,N]\,\cX[L,M]\right)v
=\lambda(a)\cdot \lambda\left(\cX[K\cap L,N\cap M]\right)v.
$$
 We have the following cases:
 
 1)
If $K\not\supset L$ or $N\not\supset M$, then by our choice of $\cX[L,M]$,
we have
$$\lambda\left(\cX[K\cap L,N\cap M]\right)=0.$$

2) Otherwise we come to
$
\lambda(a)\lambda\left(\cX[L,M]\right)v=\lambda(a)v
$.

2.1)
 If $a\in \wh \G_m[L,M]$, we get $\lambda(a)v\in V$.

2.2)
Let $a\notin \wh \G_m[L,M]$. Then  
\begin{multline}
\lambda\left(\cX[L,M]\right)\lambda(a)\lambda\left(\cX[L,M]\right)v=
\lambda(a)\Bigl\{\lambda(a^{-1}) \lambda\left(\cX[L,M]\right)\lambda(a)\Bigr\} \lambda\left(\cX[L,M]\right)v=
\\=
\lambda(a) \lambda\left(\cX[La,M a^{t-1}]\right) 
\lambda\left(\cX[L;M]\right)v=\\=
\lambda(a) \lambda\left(\cX[La\cap L,Ma^{t-1}\cap M]\right) 
v=0.
\label{eq:prod-proj}
\end{multline}
Since an idempotent $\cX[a^{-1}L\cap L,a^tM\cap M]$
is strictly smaller than  $\cX[L,M]$,
the  $\lambda(\cX[\dots])=0$.
\hfill $\square$

\sm

{\sc End of proof of Lemma \ref{l:description}.}
Thus $H_m$ is an orthogonal direct sum of
spaces $V_x$, where $x$ ranges in the homogeneous space
 $\wh \G_m[L,M]\setminus\GL(m)$, and
$\lambda(a)$ sends each $V_x$ to $V_{xa}$. This means that
the representation $\lambda$ of $\GL(m)$ 
is induced from the representation of $\wh \G_m[L,M]$ in $V$, see., e.g., 
\cite{Ser-rep}, Sect.7.1.

Operators 
$$\lambda\bigl(\cX[La^{-1},a^{t}M]\bigr)=
\lambda(a)\lambda(\cX[L,M])\lambda(a^{-1})
$$
act as orthogonal projectors to $V_{x_0 a}$. 
A projector $\lambda\bigl(\cX[K,N]\bigr)$ is identical
on $V_{x_0a}$ if and only if $\cX[K,N]\,\cX[La^{-1},Ma^{t}]=\cX[La^{-1},Ma^{t}]$
and this give us the action of the
 semigroup of idempotents.

It remains to show the representation of $\wh \G_m[L,M]$ in $V$
is irreducible.
Assume that it contains a $\wh \G_m[L,M]$-invariant subspace $W$,
then each $V_x$ contains a copy $W_x$ of $W$ and
$\oplus_x W_x$ is a $\GL(m)$-invariant subspace in the whole $H_m$.
\hfill $\square$

\begin{corollary}
	\label{cor:last-}
Let $\lambda(\frg)$ be a nonzero operator leaving $V$
invariant. Then there is $b\in \wh \G_m[L,M]$ such that
$$
\lambda(\frg)\Bigr|_V=\rho(b)\Bigr|_V.
$$	
\end{corollary}

{\sc Proof.} This operator can be represented as $\lambda(a)\lambda(\cX[N,K])$
An operator $\lambda(\cX[N,K])$ restricted to $V$ is zero 0 or 1. Let this operator be 1.
Then $\lambda(a)$ preserves $V$ only if $a\in \wh \G_m[L,M]$. In this case
we set $b=a$.
\hfill $\square$

\sm 

Keeping in  mind Lemma \ref{l:GNS} we get the following statement:

\begin{corollary}
	\label{cor:last}
An irreducible $*$-representation of the category $\cK$ 
is determined by its height $m$, a minimal idempotent $\cX[L,M]$
acting nontrivially in $H_m$ and an irreducible representation 
$\tau$ of the group $\wh \G_m[L,M]$.
\end{corollary}

We do not claim an existence of representation corresponding to given data 
of this kind.

\sm

{\bf \punct End of proof.} Let $\rho$ be an irreducible unitary
representation of $\G$ of height $m$
 in a Hilbert space $H$.
 Then we have a chain of subspaces in $H$:
 $$
  H_{m}\longrightarrow H_{m+1}
 \longrightarrow H_{m+2}\longrightarrow\dots
 $$
  Lemma \ref{l:embed} defines a chain of semigroups
 $$
\Gamma(m)\longrightarrow\Gamma(m+1)
  \longrightarrow\Gamma(m+2)\longrightarrow\dots.
 $$
 Each semigroup $\Gamma(n)$ acts in $H$ as follows:
 in $H_n$ it acts by operators $\wt\rho_{nn}(\cdot)$,
 on $H_n^\bot$ these operators are zero
 (see Lemma \ref{l:P-alpha}).  

 On the other hand, we have a chain of groups
 $$
\GL(m)\longrightarrow\GL(m+1)\longrightarrow\GL(m+2)
 \longrightarrow\dots
 $$
 acting by unitary operators, their inductive limit is the group
 $\G_\fin$.
 Each group $\GL(n)$ preserves the subspace $H_n$, on this subspace the action
 of $\GL(n)$ coincides with the action of the group $\Aut_\cK(n)=\GL(n)$.

Consider the  data listed in Corollary \ref{cor:last}. 
We regard the subspace $V=\im \wt\rho_{mm}(\cX[L,M])\subset H_m$ as
a subspace in $H$. 
Denote  the $\GL(n)$-cyclic of $V$ by $W_n$, it is a subspace in $H_n$.

\begin{lemma}
	Let $g\in \GL(n)$. 
	If $g\in \wt \G_n[L;M]$, then $\rho(g) V=V$. Otherwise,
	$\rho(g) V\bot V$.
\end{lemma}

{\sc Proof.} In the first case, we have
$$
 \wt\rho_{nn}(\cX[L,M])\, (\wt\rho_{nn}(g))^{-1}\,
\wt\rho_{nn}(\cX[Lg,Mg^{t-1}])=\wt\rho_{nn}(\cX[L,M])
$$
and therefore the image $V$ of $\wt\rho_{nn}(\cX[L,M])$
is invariant with respect to $\rho(g)$.

In the second case we repeat the line (\ref{eq:prod-proj}).
\hfill$\square$


%

\sm 

Thus the representation of $\GL(n)$ in $W_n$ is induced from
the subgroup $\wh \G_n[L,M]$. If $k>n$, then we have embeddings
$$\GL(n)\to\GL(k), \quad \wh \G_n[L,M]\to \wh \G_k[L,M]$$
 and therefore 
the map of homogeneous spaces
$$\Xi_{n,k}:\wh \G_n[L,M]\setminus\GL(n)\to \wh \G_k[L,M]\setminus\GL(k).$$
 On the other hand, we have an embedding $W_n\to W_k$ regarding
 the orthogonal decompositions of these spaces into copies of $V$,
  therefore the map $\Xi_{n,k}$
  is an embedding.

 Finally, we get a representation of $\G_\fin$ induced from the subgroup
 $\wh \G_\fin[L,M]$. By continuity, $\G$  acts regarding the same 
 orthogonal decomposition $\oplus V_{xa}$. Hence a representation
 of $\G$ is induced from closure%
 \footnote{This closure contains $\G(m)$ and we refer to Lemma \ref{l:3m}.}
  of $\wh \G_\fin[L,M]$, i.e., $\wh\G[L,M]$. 
 
 \begin{lemma}
 	The image of $\wh \G_\fin[L,M]$ in the group of operators in $V$ coincides with the image of $\wh\G_m[L,M]$.
 \end{lemma}

{\sc Proof.} Let $u\in \wh \G_n[L,M]$. Then 
$\rho(u)$ preserves $V\subset H_m$. Therefore
$$
\rho(u)\Bigr|_{V}= P_m \rho(u) P_m\Bigr|_{V}=
\wt\rho([u]_{mm})\Bigr|_{V}.
$$
By Corollary \ref{cor:last-}, this  operator has the form
$ \rho(u')\Bigr|_{V}$, where $u'\in \wh \G_m[L,M]$.
\hfill $\square$

\sm 

Thus the representation $\tau$ of $\wh \G_\fin[L;M]$ in $V$ has a finite image.
It continuous extension to $\G[L;M]$ has the same image. The kernel of the representation $\tau$ is a closed subgroup.
Since it has a finite index, it is open. 
By Proposition \ref{pr:no-finite}, $\tau$ is trivial on the subgroup $\G^\bullet[L;M]$.

\tt
\noindent
 Math. Dept., University of Vienna/c.o Pauli Institute \\
 \&Institute for Theoretical and Experimental Physics (Moscow); \\
 \&MechMath Dept., Moscow State University;\\
 \&Institute for Information Transmission Problems;\\
 URL: http://mat.univie.ac.at/$\sim$neretin/

\end{document}